\documentclass[11pt]{article}
\usepackage{graphicx}
\usepackage{amsmath,amssymb,amsthm,amsfonts}
\usepackage{amssymb}
\usepackage{appendix}
\newtheorem{thm}{Theorem}[section]

\newtheorem{lem}[thm]{Lemma}

\theoremstyle{definition}
\newtheorem{defn}{Definition}[section]
\theoremstyle{remark}

\numberwithin{equation}{section}

\DeclareMathSymbol{\C}{\mathalpha}{AMSb}{"43}

\newcommand{\eps}{\varepsilon}

\newcommand{\alp}{\alpha}

\newcommand{\R}{{\mathbb{R}}}

\newcommand{\inte}{\int_{\mathbb{R}^N}}

\newcommand{\bsub}{\begin{subequations}}
\newcommand{\esub}{\end{subequations}$\!$}

\begin{document}
\title{Concentration Behavior of Ground States for $L^2$-Critical Schr\"{o}dinger  Equation with a Spatially Decaying Nonlinearity}
\author{Yong Luo\thanks{School of Mathematics and Statistics, and Hubei Key Laboratory of
Mathematical Sciences, Central China Normal University, P.O. Box 71010, Wuhan 430079,
P. R. China.  Email: \texttt{yluo@mail.ccnu.edu.cn}. Y. Luo is partially supported by the Project funded by China Postdoctoral Science Foundation No. 2019M662680.}
\, and\, Shu Zhang\thanks{School of Mathematics and Statistics, Central China Normal University, P.O. Box 71010, Wuhan 430079,
P. R. China.  Email: \texttt{shu@mails.ccnu.edu.cn}. 
	}
}

\date{\today}

\smallbreak \maketitle

\begin{abstract}
We consider the following  time-independent nonlinear $L^2$-critical Schr\"{o}dinger equation
\[
-\Delta u(x)+V(x)u(x)-a|x|^{-b}|u|^{1+\frac{4-2b}{N}}=\mu u(x)\,\ \hbox{in}\,\ \R^N,
\]
where $\mu\in\R$, $a>0$, $N\geq 1$, $0<b<\min\{2,N\}$, and $V(x)$ is an external potential. It is shown that ground states of the above equation can be equivalently described by minimizers of the corresponding minimization problem. In this paper, we prove that there is a threshold $a^*>0$ such that minimizer exists for $0<a<a^*$ and minimizer does not exist for any $a>a^*$. However if $a=a^*$, it is proved that whether minimizer exists depends sensitively on the value of $V(0)$. Moreover, when there is no minimizer at threshold $a^*$, we give a detailed
concentration behavior of minimizers as $a\nearrow a^*$, based on which we finally prove that there is a unique minimizer as $a\nearrow a^*$.

%
%
%
\end{abstract}

\vskip 0.2truein


\vskip 0.2truein

\section{Introduction}
In this paper, we study the following time-independent nonlinear focusing Schr\"{o}dinger  equation:
\begin{equation}\label{F}
-\Delta u(x)+V(x)u(x)-a u^{1+2\beta^2}(x)|x|^{-b}=\mu u(x)\ \ \mbox{in}\ \ \R^{N},
\end{equation}
where $N\geq1$, $a,\,\beta>0$, $0<b<\min\{2,N\}$, $\mu\in\R$ is the chemical potential, and $V(x)$ is an external potential. Equation \eqref{F} arises in various physical context such as nonlinear optics, plasma physics and Bose-Einstein condensates (BECs) (cf. \cite{BGC,Liu,Za}), where the function $a|u|^{2\beta^2}|x|^{-b}$ is a potential which can either stand for corrections to the nonlinear power-law response or for some inhomogeneities
in the medium. Due to the singularity of $|x|^{-b}$ at $x=0$, much attention have been attracted to existence, asymptotic behavior and uniqueness of solutions for \eqref{F} in the last two decades, see \cite{AD,GF1,FG,SC,Tol} and the references therein.

As illustrated in Theorem \ref{A.3} in the Appendix, ground states of \eqref{F} can be equivalently described by minimizers of
\begin{equation}\label{e(a)}
e(a)=\inf \limits_{u\in\mathcal{M}}E_{a}(u),
\end{equation}
where the energy functional $E_{a}(u)$ is defined by
\begin{equation}\label{Eau}
E_{a}(u):=\int_{\R^{N}}(|\nabla u|^{2}+V(x)|u(x)|^{2})dx-\frac{a}{1+\beta^2}\int_{\R^{N}}\frac{|u(x)|^{2+2\beta^2}}{|x|^{b}}dx,\,\ u\in\mathcal{H},
\end{equation}
and
\begin{equation}\label{M}
\mathcal{M}:=\big\{u\in\mathcal{H}:\int_{\R^{N}}|u(x)|^{2}dx=1\big\}.
\end{equation}
Here the space $\mathcal{H}$ is defined as
\begin{equation}
\mathcal{H}:=\Big\{u\in H^{1}(\R^{N}):\int_{\R^{N}}V(x)|u(x)|^{2}dx<\infty\Big\}
\end{equation}
with the associated norm $\|u\|_{\mathcal{H}}=\Big\{\displaystyle\int_{\R^N}\Big(|\nabla u(x)|^{2}+[1+V(x)]|u(x)|^{2}\Big)dx\Big\}^{\frac{1}{2}}$. To discuss equivalently ground states of \eqref{F}, throughout the whole paper we shall therefore focus on investigating \eqref{e(a)}, instead of \eqref{F}.


In the sequel, we always assume that $V(x)$ satisfies
\begin{itemize}
\item[\rm($V_{1})$.] $0\leq V(x)\in {L^{\infty}_{loc}} (\R^{N})$  and $ \lim_{|x|\rightarrow\infty}V(x)=\infty$.
\end{itemize}
When $b=0$, (\ref{e(a)}) is a homogeneous constraint variational problem, for which there are many existing results of (\ref{e(a)}) (cf. \cite{GTS,MP,GR,GW,GZ,ZJ}). To be more precise, when $\beta^2>\frac{2}{N}$, (\ref{e(a)}) is the so-called $L^2$-supercritical problem. One can take a suitable test function to derive that $e(a)=-\infty$, and thus  (\ref{e(a)}) admits no minimizers for any $a>0$ under this case. However, if $0<\beta^2<\frac{2}{N}$, (\ref{e(a)}) is in the $L^2$-subcritical case, and we can prove \eqref{e(a)} always
admits minimizers for any $a>0$ by using energy estimates and the compactness embedding in Lemma \ref{lem2.1} below, see also \cite{SX}. As for the $L^2$-critical case where $\beta^2=\frac{2}{N}$, the existence and nonexistence of minimizers for \eqref{e(a)} and various quantitative properties of \eqref{e(a)} have been investigated by Guo and his co-authors in recent years, see \cite{GLW,GR,GW,GZ}.

When $b\neq0$,  (\ref{e(a)}) is the inhomogeneous constraint  variational problem which contains the nonlinear term $K(x)|u(x)|^{2+2\beta^2}$ with $K(x)=|x|^{-b}$. Note that similar inhomogeneous problems were analyzed recently in \cite{DGL1,DGL2}, which however focused mainly on the case where $K(x)\in L^\infty(\R^N)$ does not have any singularity. On the other hand, \cite{FG} studied the existence, instability of standing waves of \eqref{F} in the $L^2$-critical case $\beta^2=\frac{2-b}{N}$ and $V(x)\equiv 1$. Inspired by the above works, in this paper we shall mainly studied minimizers (ground states) of \eqref{e(a)} in the critical case $\beta^2=\frac{2-b}{N}$, where $0<b<\min\{2,N\}$. For a general class of $V(x)$ satisfying the assumptions $(V_1)$, we shall give the existence and nonexistence, blow up behavior, and local uniqueness of minimizers for \eqref{e(a)}.

We first introduce the following associated limit equation
\begin{equation}\label{Eq}
-\triangle u+u-|x|^{-b}|u|^{2\beta^2}u=0 \ \mbox{in}\ \ \R^{N}.
\end{equation}
Note that \eqref{Eq} admits a unique positive radial solution, see the proof in \cite{GF1} for $N\geq 3$, \cite{GF2} for $N=2$, and \cite{Tol} for $N=1$. Thus we always denote this unique positive solution by $Q(x)$ in this paper. Moreover,
we cite from \cite[Corollary 2.3]{FG} the following Gagliardo-Nirenberg inequality,
\begin{equation}\label{GN}
\frac{a^*}{1+\beta^2}\int_{\R^{N}}|x|^{-b}|u(x)|^{2+2\beta^2}dx\leq \|\nabla u(x)\|^{2}_{2}\|u(x)\|^{2\beta^2}_{2} \ \mbox{in}\ \ H^{1}(\R^{N}),
\end{equation}
where the constant $a^*>0$ is given by
\begin{equation}\label{a8}
a^*:=\|Q\|^{2\beta^2}_{2},
\end{equation}
and $``="$ holds in \eqref{GN} if and only if $u(x)=mn^{\frac{N}{2}} Q(nx)$ ($m,n\neq 0$ are arbitrary).
Recall also from \cite[Theorem 2.2]{SC} that $Q$ admits the following exponential decay,
\begin{equation}\label{zs}
|Q(x)|,\quad |\nabla Q(x)|\leq Ce^{-|x|}\ \ \mbox{as}\ \ |x|\rightarrow\infty.
\end{equation}
Moreover, one can derive from (\ref{Eq}) and (\ref{GN}) that $Q$ satisfies
\begin{equation}\label{DS}
\int_{\R^{N}}|\nabla Q|^{2}dx=\frac{1}{\beta^2}\int_{\R^{N}}|Q|^{2}dx=\frac{1}{1+\beta^2}\int_{\R^{N}}|x|^{-b}|Q|^{2+2\beta^2}dx.
\end{equation}

By applying the Gagliardo-Nirenberg inequality (\ref{GN}), we can obtain the following theorem concerning the existence and non-existence of the minimizers for (\ref{e(a)}).
\vskip 0.05truein

\begin{thm}\label{th1}
Let $N\geq 1$, $0<b<\min\{2,N\}$, $\beta^2=\frac{2-b}{N}$, and assume $V(x)$ satisfies $(V_{1})$, then
\begin{enumerate}
\item If $0\leq a< a^{\ast}=\|Q\|^{2\beta^2}_{2}$, there exists at least one minimizer for (\ref{e(a)}).
\item If $a>a^{\ast}$, there is no minimizer for (\ref{e(a)}).
\end{enumerate}
Moreover, $e(a)>0$ for $a<a^{\ast}$, and $e(a)=-\infty$ for $a>a^{\ast}$.
\end{thm}
In view of Theorem \ref{th1}, a natural question arises whether there exist minimizers for $e(a)$ when $a=a^*$.
In contrast to the homogeneous $b=0$, 
our following results show that $e(a)$ may admit minimizers at the threshold $a^*$ for some trapping potentials $V(x)$ satisfying $(V_1)$.
Before state our results, we first give some useful energy estimates for $e(a^*)$ under the assumption of $(V_1)$. Using the Gagliardo-Nirenberg inequality, it is easy to obtain that $e(a^*)\geq \inf_{x\in\R^{N}}V(x)$. On other hand, taking the same test function as in \eqref{st} below
and letting $\tau\to\infty$, one can obtain that $e(a^*)\leq E(u_\tau)\leq V(0)$. Therefore, we conclude from above that 
$$\inf_{x\in\R^{N}}V(x)\leq e(a^*)\leq V(0),$$ 
and thus our second result can be summarized as follows:
\begin{thm}\label{th2}
Under the assumptions of Theorem \ref{th1} and assume $a=a^*$, then we have
\begin{enumerate}
\item If $\inf_{x\in\R^{N}}V(x)\leq e(a^*)<V(0)$, then there exists at least one minimizer for $e(a^{\ast})$.
\item If  $\inf_{x\in\R^{N}}V(x)= e(a^*)=V(0)$, then there is no minimizer for $e(a^{\ast})$. Moreover, $\lim_{a\nearrow a^{\ast}}e(a)=e(a^{\ast})=V(0)$.
\end{enumerate}
\end{thm}
Theorem \ref{th2} tells us that the existence and non-existence of minimizers for \eqref{e(a)} depend sensitively on the value of $e(a^*)$ and $V(0)$. Heuristically speaking, $V(0)$ can be regarded as a threshold of the energy $e(a^*)$ such that \eqref{e(a)} admits minimizer if $e(a^*)$ is less than $V(0)$,
while \eqref{e(a)} does
not admit any minimizer if $\inf_{x\in\R^{N}}V(x)= e(a^*)=V(0)$.


In view of Theorem \ref{th2} (2), we obtain that $u_a$ blow up in the sense that $\inte |\nabla u_a|^2dx\to\infty$ as $a\nearrow a^*$. Therefore, we next focus on the limit behavior of minimizers $u_a$ as $a\nearrow a^*$. Note that $E_{a}(u_a)= E_{a}(|u_a|)$ \cite{LE}, therefore $|u_a|$ is also a minimizer of $e(a)$. By the strong maximum principle, we can further derive from (\ref{F}) that $|u_a|>0$ holds in $\R^N$. Therefore, without loss of generality, we can restrict the minimizers of $e(a)$ to positive functions.


To state our results, we need some additional assumptions of $V(x)$, for which we define
\begin{defn}\label{def1}
The function $h(x)\geq0$ in $\R^{N}$ is homogeneous of degree $l\in \R^{+}$(with respect to the origin), if there exists some $l>0$ such that
\begin{equation}\label{v1}
h(tx)=t^{l}h(x)\ \ \mbox{in}\ \ \R^{N}\ \ \mbox{ for any}\ \ t>0.
\end{equation}
\end{defn}
One can obtain from above definition that the homogeneous function $h(x)\in C(\R^{N})$ of degree $l>0$ satisfies
$0\leq h(x)\leq C |x|^{l}$ in $\R^{N}$, where $C=\max\limits_{x\in \partial B_{1}(0)}h(x)$. Moreover, if the homogenous function $h(x)$ satisfies $\lim_{|x|\rightarrow\infty}h(x)=+\infty$, then $x=0$ is the unique minimum point of $h(x)$.

Following the above definition and the assumption $(V_1)$, we next assume that 
\begin{enumerate}
\item [\rm($V_2$).]
$V(x)\in C^1(\R^2)$ satisfies $\{x\in\R^N:\,V(x)=0\}=\{0\}$, 
\begin{equation}\label{s.3}
|V(x)|\le C e^{\gamma|x|},\,\ |\nabla V(x)|\le C e^{\gamma|x|}\quad\hbox{for some $\gamma>0$ as $|x|\to\infty$},
\end{equation}
and for $j=1,2,\cdots,N,$
\begin{equation}\label{s.4}
V(x)=h(x)+o(|x|^l),\ \,\frac{\partial V(x)}{\partial{x_j}}=
\frac{\partial h(x)}{\partial{x_j}}+o(|x|^{l-1})\,\ \hbox{as}\,\  |x|\to 0,
\end{equation}
where $0\le h(x)\in C^{1}(\R^N)$ is a homogeneous function of degree $l>0$ and satisfies $\lim_{|x|\to \infty}h(x)=+\infty$.
\end{enumerate}
For simplicity, we also denote
\begin{equation}\label{d2}
\lambda:=\Big(\frac{a^{\ast}l}{2 \beta^{l}\|Q\|^{2}_{2}}\int_{\R^{N}}h(x)Q^{2}(x)dx\Big)^{\frac{1}{l+2}}.
\end{equation}

Combined the assumption $(V_1)$ with $\{x\in\R^N:\,V(x)=0\}=\{0\}$ yields that $0=V(0)=\inf_{x\in\R^{N}}V(x)$, and hence the assumptions of Theorem \ref{th2} (2) hold true. 
We give the detailed limit behavior of minimizers $u_a$ as follows:
\begin{thm}\label{th3}
Under the assumptions of Theorem \ref{th1} and suppose that $V(x)$ satisfies $(V_2)$. Then for any positive minimizer $u_{a}$ of (\ref{e(a)}), we have
\begin{equation}\label{bp}
\lim\limits_{a\nearrow a^*} \frac{(a^{\ast}-a)^{\frac{N}{2l+4}}\|Q\|_{2}}{(\lambda\beta)^{\frac{N}{2}}}u_{a}\Big(\frac{\big(a^{\ast}-a\big)^{\frac{1}{l+2}}}{\lambda\beta}x\Big)=Q(x)\quad\hbox{}
\end{equation}
strongly in $L^{\infty}(\R^{N})\cap H^1(\R^N)$, where $\lambda$ is defined by (\ref{d2}).
\end{thm}
Theorem \ref{th3} shows that minimizers $u_a$ of $e(a)$ must concentrate at $0$ when $V(0)=\inf_{x\in\R^{N}} V(x)=0$. Theorem \ref{th3} is proved by employing the refined energy estimates and blow up analysis. Compared with those homogeneous problems, we should point out that our analysis on inhomogeneous problem \eqref{e(a)} are more complicated. Moreover, there is some essential difficulties on obtain the $L^\infty$-uniform convergence of \eqref{bp}
due to the singularity of $|x|^{-b}$. 

Motivated by  \cite{Cao,Deng,GLW}, we finally investigate the uniqueness of positive minimizers for $e(a)$ as $a\nearrow a^*$.
\begin{thm}\label{unique}
Under the assumptions of Theorem \ref{th3} and suppose $N\geq 3$, then there exists a unique positive minimizer for $e(a)$ as $a\nearrow a^*$.
\end{thm}

We remark that the restriction $N\geq 3$ can be removed in Theorem \ref{unique} if the non-degeneracy property \eqref{2:linearized} below still holds true for any dimension $N\geq 1$.
Similar to those existing works in \cite{Cao,Deng,GLW},
Theorem \ref{unique} is proved by establishing local Pohozaev identity.  However, the calculations included in the proof is more involved due to the inhomogeneous nonlinear term. Besides, we mention that the standard elliptic regularity theory should be used carefully when concerning the singular term near  the origin.

This paper is organized as follows. In Section 2 we shall prove Theorems \ref{th1} and \ref{th2} on
the existence and non-existence of minimizers for (\ref{e(a)}). Section 3 is then devoted to the proof of Theorem \ref{th3} on the concentration behavior of minimizes for (\ref{e(a)}). By establishing local Pohozaev identity, we finally prove Theorem \ref{unique} on the local uniqueness of positive minimizers in Section 4.
\section{Existence of minimizers}
This section is devoted to prove Theorems 1.1 and 1.2 with regard to the existence and nonexistence of minimizers for the problem (\ref{e(a)}). Motivated by \cite{GR}, we first introduce the following compactness lemma \cite[Theorem 2.1]{TZ}.
\begin{lem}\label{lem2.1}
	Assume $V(x)\in L^{\infty}_{loc}(\R^{N})$ and $\lim_{|x|\rightarrow\infty}V(x)=\infty$, then the embedding $\mathcal{H}\hookrightarrow L^{q}(\R^{N})$ is compact for $2\leq q <2^*$, where $2^*=\frac{2N}{N-2}$ if $N\geq 3$ and $2^*=+\infty$ if $N=1,2$.
\end{lem}
This compactness property, together with the Gagliardo-Nirenberg inequality (\ref{GN}), allows us to prove Theorem \ref{th1}.

\vskip 0.05truein
\noindent{\bf Proof of Theorem \ref{th1}.}

(1). Firstly, we prove that $e(a)$ has at least one minimizer for all $0\leq a< a^{\ast}=\|Q\|^{2\beta^2}_{2}$.
Indeed, for any fixed  $0\leq a<a^{\ast}$ and $u\in \mathcal{M}$, we derive from the assumption $(V_1)$ and Gagliardo-Nirenberg inequality (\ref{GN}) that
\begin{equation}\label{Eu}
\begin{split}
E_{a}(u)&=\int_{\R^{N}}(|\nabla u(x)|^{2}+V(x)|u(x)|^{2})dx-\frac{a}{1+\beta^2}\int_{\R^{N}}|u(x)|^{2+2\beta^2}|x|^{-b}dx\\
& \geq(1-\frac{a}{a^*})\int_{\R^{N}}|\nabla u(x)|^{2}dx+\int_{\R^{N}}V(x)|u(x)|^{2}dx\\
&\geq(1-\frac{a}{a^*})\int_{\R^{N}}|\nabla u(x)|^{2}dx\geq 0.
\end{split}
\end{equation}
This implies that $e(a)\geq0$ for $0<a<a^*$.

Then we can choose a minimizing sequence $\{u_{n}\}\subset\mathcal{H}$ satisfying $\|u_{n}\|_{2}=1$ and $\lim\limits_{n\rightarrow\infty}E_{a}(u_{n})=e(a)$. Due to (\ref{Eu}), we obtain that $\int_{\R^{N}}|\nabla u_{n}(x)|^{2}dx$ and $\int_{\R^{N}}V(x)|u_{n}|^{2}dx$ are bounded uniformly for $n$. Following Lemma 2.1, we may assume passing if necessary to a subsequence,
$$u_{n}\rightharpoonup u\,\ \hbox{weakly in} \,\ \mathcal{H},\,\  u_{n}\rightarrow u\,\ \hbox{strongly in}\,\ L^{q}(\R^{N}),\,\ 2\leq q<2^* $$
for some $u\in \mathcal{H}$.
Thus we conclude that $\int_{\R^{N}}|u(x)|^{2}dx=1$.
For any $s\in (\frac{N}{N-b},\infty)$, note that any
function $f\in L^1(\R^N)\cap L^s(\R^N)$ satisfies
\begin{equation}\label{2.4}
\begin{aligned}
\inte f|x|^{-b}dx
&\leq \|f\|_{L^{s}(B_R(0))}\big\||x|^{-b}\big\|_{L^{\bar s}(B_R(0))}+\frac{1}{R^b}\|f\|_{L^{1}(\R^N)}\\
&\leq C(R)\|f\|_{L^{s}(B_R(0))}+\frac{1}{R^b}\|f\|_{L^{1}(\R^N)},
\end{aligned}
\end{equation}
where $\bar s=\frac{s}{s-1}$ satisfies $1<\bar s<\frac{N}{b}$. Therefore, taking $f_n=|u_{n}(x)|^{2+2\beta^2}-|u(x)|^{2+2\beta^2}$, we have $f_n\to 0$ in $\in L^1(\R^N)\cap L^s(\R^N)$ for $s\in(\frac{N}{N-b},\frac{N^2}{(N-2)(N+2-b)})$ if $n\geq 3$, $s\in (\frac{N}{N-b},\infty)$ if $N=1,2$, and hence
\begin{equation}\label{2.1}
\lim_{n\rightarrow\infty}\int_{\R^{N}}|x|^{-b}|u_{n}(x)|^{2+2\beta^2}dx=\int_{\R^{N}}|x|^{-b}|u(x)|^{2+2\beta^2}dx.
\end{equation}
Consequently, by weak lower semicontinuity, we deduce that $e(a)=\lim_{n\to\infty}E_{a}(u_n)\geq E_{a}(u)\geq e(a)$. Therefore, $u$ must be a minimizer of $e(a)$, thus Theorem \ref{th1} (1) is proved and $e(a)=E_a(u)>0$.

(2). We next show that there is no minimizer for $e(a)$ as $a>a^{\ast}$.
Let $\varphi\in C_{0}^{\infty}(\R^{N})$ be a nonnegative cut-off function  satisfying $\varphi(x)=1$ if $|x|\leq1$, and $\varphi(x)=0$ if $|x|>2$. Set for all $\tau>0$,
\begin{equation}\label{st}
u_{\tau}(x)=A_{\tau}\frac{\tau^{\frac{N}{2}}}{\|Q\|_{2}}\varphi(x)Q(\tau x),
\end{equation}
where $A_{\tau}>0$ is chosen so that $\int_{\R^{N}}|u_{\tau}(x)|^{2}dx=1$. By scaling, we deduce from the exponential decay of $Q$ in  (\ref{zs}) that
\begin{equation}\label{A}
\frac{1}{A^{2}_{\tau}}=\frac{1}{\|Q\|^{2}_{2}}\int_{\R^{N}}Q^{2}(x)\varphi^{2}(\frac{x}{\tau })=1+o\big((\tau)^{-\infty}\big) \ \ \mbox{as}\ \ \tau\rightarrow\infty.
\end{equation}
Here we use the notation $f(t)=o(t^{-\infty})$ for a function $f$ satisfying $\lim\limits_{t\rightarrow\infty}|f(t)|t^{s}=0$ for all $s>0$.

Following (\ref{zs}) and (\ref{DS}), we also obtain that
\begin{equation}\label{Gj}
\begin{split}
&\int_{\R^{N}}|\nabla u_{\tau}(x)|^{2}dx-\frac{a}{1+\beta^2}\int_{\R^{N}}|x|^{-b}u_{\tau}^{2+2\beta^2}dx
\\&=\frac{\tau^{2}A^{2}_{\tau}}{\|Q\|^{2}_{2}}\int_{\R^{N}}|\nabla Q|^{2}dx-\frac{a}{1+\beta^2}\frac{\tau^{2}A^{2+2\beta^2}_{\tau}}{\|Q\|^{2+2\beta^2}_{2}}\int_{\R^{N}}|x|^{-b}Q(x)^{2+2\beta^2}dx+o(\tau^{-\infty})
\\&=\frac{\tau^{2}A^{2}_{\tau}}{\|Q\|^{2}_{2}}\int_{\R^{N}}|\nabla Q|^{2}dx-\frac{a}{1+\beta^2}\frac{\tau^{2}A^{2+2\beta^2}_{\tau}}{\|Q\|^{2+2\beta^2}_{2}}(1+\beta^2)\int_{\R^{N}}|\nabla Q|^{2}dx+o(\tau^{-\infty})
\\&=\frac{\tau^{2}}{\|Q\|^{2}_{2}}(1-\frac{a}{a^*})\int_{\R^{N}}|\nabla Q|^{2}dx+o(\tau^{-\infty})\ \ \mbox{as}\ \ \ \tau\rightarrow\infty.
\end{split}
\end{equation}
On the other hand, since the function $V(x)\varphi^{2}(x)$ is bounded and has compact support, it follows from \cite{LE} that
\begin{equation}\label{jx}
\lim\limits_{\tau\rightarrow\infty}\int_{\R^{N}}V(x)|u_{\tau}(x)|^{2}dx=V(0).
\end{equation}
Hence, we have
\begin{equation}\label{zh}
\begin{split}
e(a)&\leq E_{a}(u_{\tau}(x))\\
&=\int_{\R^{N}}|\nabla u_{\tau}(x)|^{2}dx-\frac{a}{1+\beta^2}\int_{\R^{N}}|x|^{-b}u^{2+2\beta^2}_{\tau}(x)dx+\int_{\R^{N}}V(x)|u_{\tau}(x)|^{2}dx\\
&=\frac{\tau^{2}}{\|Q\|^{2}_{2}}(1-\frac{a}{a^*})\int_{\R^{N}}|\nabla Q|^{2}dx+\int_{\R^{N}}V(x)|u_{\tau}(x)|^{2}dx+o(\tau^{-\infty})\\
&\quad\longrightarrow-\infty\quad\quad\hbox{as}\quad \tau\to\infty.
\end{split}
\end{equation}
Therefore, if $a>a^{\ast}$, $e(a)=-\infty$, which implies the non-existence of minimizers as $a>a^{\ast}$.
\qed

\noindent\textbf{Proof of Theorem \ref{th2}.}
(1). Following the Gagliardo-inequality \eqref{GN}, one can check that $e(a^*)\geq 0$ is bounded from below, and thus there exists a minimizing sequence $\{u_{n}\}$ of $e(a^{\ast})$ such that $e(a^*)=\lim_{n\to\infty}E_a(u_n)$.
To prove Theorem \ref{th2} (1), by the compact embedding in Lemma 2.1, to prove Theorem \ref{th2}, it is enough to prove that $\{u_{n}\}$ is bounded uniformly in $\mathcal{H}$. On the contrary, we assume that $\{u_{n}\}$ is unbounded in $\mathcal{H}$, then there exists a subsequence of $\{u_{n}\}$, still denoted by $\{u_{n}\}$, such that $\|u_{n}\|_{\mathcal{H}}\rightarrow\infty$ as $n\rightarrow\infty$. Moreover, by the Gagliardo-Nirenberg inequality, we have
\begin{equation}\label{V}
\int_{\R^{N}}V(x)|u_{n}(x)|^{2}dx\leq E_{a^{\ast}}(u_{n})\leq e(a^{\ast})+1.
\end{equation}
Since $\{u_{n}\}$ is unbounded in $\mathcal{H}$, then we have
\begin{equation}\label{U}
\int_{\R^{N}}|\nabla u_{n}(x)|^{2}dx\rightarrow\infty.
\end{equation}
Let
\begin{equation}\label{wq}
\varepsilon^{-2}_{n}:=\int_{\R^{N}}|\nabla u_{n}(x)|^{2}dx,
\end{equation}
then by (\ref{U}), we have that $\varepsilon_{n}\rightarrow 0$ as $n\rightarrow\infty$. Now we defined the $L^{2}(\R^{N})-$normalized function
\begin{equation}\label{omg}
w_{n}(x):=\varepsilon^{\frac{N}{2}}_{n}u_{n}(\varepsilon_{n}x),
\end{equation}
which satisfies
\begin{equation}\label{2gj}
\int_{\R^{N}}|\nabla w_{n}(x)|^{2}dx=\int_{\R^{N}}| w_{n}(x)|^{2}dx=1.
\end{equation}
Therefore $\{w_n\}$ is bounded uniformly in $H^1(\R^N)$ for $n$, and Sobolev embedding theorem implies that passing if necessary to a subsequence
$$w_{n}\rightharpoonup w_0\,\ \hbox{weakly in} \,\ H^1(\R^N),\,\  w_{n}\rightarrow w_0\,\ \hbox{strongly in}\,\ L_{loc}^{q}(\R^{N}),\,\ 2\leq q<2^*. $$

Next, we claim that
\begin{equation}\label{2.2}
\lim_{n\to\infty}w_n=\frac{\beta^{\frac{N}{2}}Q(\beta x)}{\|Q\|_2}\quad\hbox{strongly in $H^1(\R^N)$}.
\end{equation}
Actually, by the definition of $w_n$, we obtain from Gagliardo-Nirenberg inequality (\ref{GN}) that
\begin{equation}\label{wgx}
\begin{split}
C&\geq e(a^*)=\lim_{n\rightarrow\infty}E_{a^*}(u_{n})
\\&=\lim_{n\rightarrow\infty}\Big\{\int_{\R^{N}}\Big(|\nabla u_{n}|^{2}-\frac{a^*}{1+\beta^2}|x|^{-b}u_{n}^{2+2\beta^2}+V(x)|u_{n}(x)|^{2}\Big)dx\Big\}
\\
&=\lim_{n\rightarrow\infty}\Big\{\frac{1}{\varepsilon^{2}_{n}}\int_{\R^{N}}\Big[|\nabla w_{n}|^{2}-\frac{a^{\ast}}{1+\beta^2}|x|^{-b}w_{n}^{2+2\beta^2}\Big]dx+\int_{\R^{N}}V(\varepsilon_{n}x)|w_{n}(x)|^{2}dx\Big\}\\
&\geq 0.
\end{split}
\end{equation}
Therefore, we derive from \eqref{2gj}, the fact $\eps_n\to 0$ as $n\to\infty$ and above that
\begin{equation}\label{wnp}
\lim_{n\to\infty}\frac{a^{\ast}}{1+\beta^2}\int_{\R^{N}}|x|^{-b}w_{n}^{2+2\beta^2}(x)dx=\lim_{n\to\infty}\int_{\R^{N}}|\nabla w_{n}(x)|^{2}dx=1.
\end{equation}
On the other hand, similar to the proof of \eqref{2.1}, we have
\begin{equation}\label{lpsl}
\frac{a^{\ast}}{1+\beta^2}\int_{\R^{N}}|x|^{-b}w_{0}^{2+2\beta^2}(x)dx=\lim_{n\to\infty}\frac{a^{\ast}}{1+\beta^2}\int_{\R^{N}}|x|^{-b}w_{n}^{2+2\beta^2}(x)dx=1,
\end{equation}
which implies $w_{0}\not\equiv0$. Moreover, by the Gagliardo-Nirenberg inequality (\ref{GN}) and weak lower semicontinuity, we have
\begin{equation}\label{wgx1}
\begin{split}
0&= \lim_{n\rightarrow\infty}\varepsilon_n^2e(a^*)
\\
&\geq\lim_{n\rightarrow\infty}\int_{\R^{N}}\Big[|\nabla w_{n}|^{2}-\frac{a^{\ast}}{1+\beta^2}|x|^{-b}w_{n}^{2+2\beta^2}\Big]dx
\\&\geq\int_{\R^{N}}|\nabla w_{0}(x)|^{2}dx-\frac{a^{\ast}}{1+\beta^2}\int_{\R^{N}}|x|^{-b}w_{0}^{2+2\beta^2}(x)dx
\\&\geq\Big[1-\Big(\int_{\R^{N}}w^2_{0}(x)dx\Big)^{\beta^2}\Big]\int_{\R^{N}}|\nabla w_{0}(x)|^{2}dx.
\end{split}
\end{equation}
Since $\|w_0\|_{L^2(\R^N)}\leq 1$, the above inequality further implies that
\begin{equation}\label{2.3}
\int_{\R^{N}}|w_{0}(x)|^{2}dx=1,\,\ \int_{\R^{N}}|\nabla w_{0}(x)|^{2}dx=\frac{2a^{\ast}}{1+\beta^2}\int_{\R^{N}}|x|^{-b}w_{0}^{2+2\beta^2}(x)dx=1.
\end{equation}
Combing the above facts, we obtain that Gagliardo-Nirenberg inequality is achieved by $w_0$ and hence $w_{0}=\frac{\beta^{\frac{N}{2}}Q(\beta x)}{\|Q\|_2}$ in view of \eqref{2.3} and \eqref{DS}.  By the norm preservation, we finally conclude that $w_{n}$ converges to $w_{0}=\frac{\beta^{\frac{N}{2}}Q(\beta x)}{\|Q\|_2}$ strongly in $H^{1}(\R^{N})$ and \eqref{2.2} is thus proved.

Based on \eqref{2.2}, by Fatou's lemma, we obtain
\begin{equation}\label{md1}
e(a^{\ast})\geq\lim_{n\rightarrow\infty}\int_{\R^{N}}V(\varepsilon_{n} x)|w_{n}(x)|^{2}dx\geq \int_{\R^{N}}V(0)|w_{0}(x)|^{2}dx=V(0),
\end{equation}
which contradicts the assumption  $e(a^{\ast})<V(0)$. Hence the minimizing sequence $\{u_{n}\}$ of $e(a^*)$ is bounded uniformly in $\mathcal{H}$ and Theorem \ref{th2} (1) is proved.

(2). We next consider the case $a=a^{\ast}$ and $V(0)=0$. Note that $e(a^*)\geq 0$. On the other hand, we can use the same test function as that of (\ref{st}) to derive that $e(a^{\ast})\leq V(0)=0$,  and thus $e(a^{\ast})=0$. Now suppose that there exists a minimizer $u\in\mathcal{M}$ at $a=a^{\ast}$. Then we have
\begin{equation}\label{1}
\int_{\R^{N}}V(x)|u(x)|^{2}dx=0=\inf\limits_{x\in\R^{N}}V(x),
\end{equation}
as well as
\begin{equation}\label{2}
\begin{split}
\int_{\R^{N}}|\nabla u(x)|^{2}dx&=\frac{a^{\ast}}{1+\beta^2}\int_{\R^{N}}|x|^{-b}|u(x)|^{2+2\beta^2}dx.
\end{split}
\end{equation}
This is a contradiction, since from the first equality $u$ must have compact support, while from the last equality it must equal to $\frac{n^{\frac{N}{2}}Q(n x)}{\|Q\|_2}>0$, where $n\in\R$ is a constant. This gives the nonexistence of minimizers for $e(a^*)$ at $a=a^*$. Moreover, taking the same test function as in (\ref{st}), we have $V(0)=\inf_{x\in\R^{N}}V(x)\leq\lim_{a\nearrow a^{\ast}}e(a)\leq V(0)$ and hence $\lim_{a\nearrow a^{\ast}}e(a)=V(0)=e(a^*)$. This completes the proof of Theorem \ref{th2}.\qed

\section{Mass concentration as $a\nearrow a^{\ast}$}
In this section, we shall prove Theorem \ref{th3}, which focuses on the concentration behavior of positive minimizers of $e(a)$ as $a\nearrow a^{\ast}$. Let $u_{a}$ be a positive minimizer of $e(a)$, by variational theory,  $u_{a}$ satisfies the following Euler-Lagrange equation
\begin{equation}\label{EL}
-\Delta u_{a}(x)+V(x)u_{a}(x)=\mu_a u_{a}(x)+a u_{a}^{1+2\beta^2}(x)|x|^{-b}\quad \mbox{in}\ \ \R^{N},
\end{equation}
where $\mu_{a}\in\R$ is a suitable Lagrange multiplier $\mu_a\in\R$ satisfying
\begin{equation}\label{2.8}
\mu_a=e(a)-\frac{a\beta^2}{1+\beta^2}\inte |x|^{-b}|u_a|^{2+2\beta^2}dx.
\end{equation}

We first establish the following lemma.

\begin{lem}\label{yl1}
Under the assumptions of Theorem \ref{th3}, and let $u_{a}$ be a positive minimizer of $e(a)$. Define
\begin{equation}\label{yp}
\varepsilon_{a}:=\Big(\int_{\R^{N}}|\nabla u_{a}|^{2}dx\Big)^{-\frac{1}{2}},
\end{equation}
and
\begin{equation}\label{xpa}
w_{a}(x)=\varepsilon^{\frac{N}{2}}_{a}u_{a}(\varepsilon_{a}x),
\end{equation}
then we have
\begin{enumerate}
\item $\varepsilon_{a}\rightarrow 0$   as  $a\nearrow a^{\ast}$.
\item    $w_a\rightarrow\frac{\beta^{\frac{N}{2}}Q(\beta x)}{\|Q\|_2}$ strongly in $ H^1(\R^N)\cap L^{\infty}(\R^{N})$ as $a\nearrow a^{\ast}$.
\item There exists a positive constant $C>0$ independent of $a\in(a,a^*)$ such that
\begin{equation}\label{zsll}
|w_a(x)|,\,\ |\nabla w_a(x)|\leq Ce^{-\frac{\beta}{4}|x|}\ \ \mbox{in}\ \ \R^N.
\end{equation}
		%
\end{enumerate}
\end{lem}
\noindent\textbf{Proof.}  1. On the contrary, suppose that there exists a sequence $\{a_{k}\}$ satisfies $a_{k}\nearrow a^{\ast}$ as $k\rightarrow\infty$ such that the sequence $\{u_{a_{k}}\}$ is bounded uniformly in $H^1(\R^N)$. By Gagliardo-Nirenberg inequality, we obtain that $\inte V(x)|u_n|^2dx\leq e(a)\leq C$ uniformly for $n$ and any $a\in(0,a^*)$ and hence $\{u_{a_{k}}\}$ is also bounded uniformly in $\mathcal{H}$.
Following Lemma \ref{lem2.1},  there exists a subsequence, still denoted by $\{a_{k}\}$ , of $\{a_{k}\}$ and $u_{0}\in\mathcal{H}$ such that
$${u_{a_{k}}}\rightharpoonup u_{0}\,\ \hbox{weakly in}\,\ \mathcal{H}\,\ \hbox{and}\,\   {u_{a_{k}}}\rightarrow u_{0}\,\ \hbox{strongly in}\,\  L^{q}(\R^{N})\,(2\leq q <2^*).$$
Combining the above convergence with \eqref{2.4}  gives that
\begin{equation*}
0=e(a^{\ast})\leq E_{a^{\ast}}(u_{0})\leq \lim\limits_{k\rightarrow\infty}E_{a_{k}}(u_{a_{k}})=\lim\limits_{k\rightarrow\infty}e(a_{k})=0,
\end{equation*}
and thus $u_{0}$ is a minimizer of $e(a^{\ast})$, which however contradicts with Theorem \ref{th2} (2). Therefore we concluded that $\varepsilon_{a}\rightarrow 0$ as $a\nearrow a^{\ast}$.

2. Note from Theorem \ref{th2} (2) that $0=e(a^*)=\lim_{a\nearrow a^*}e(a)$, thus $\{u_a\}$ is also a minimizing sequence of $e(a^*)$. Moreover, we obtain from \eqref{yp} and \eqref{xpa} that
\[
\inte |\nabla w_a|^2dx=\inte |w_a|^2dx=1\,\ \hbox{for any $a\nearrow a^*$.}
\]
Hence $\{w_a\}$ is bounded uniformly in $H^1(\R^N)$ for $a$, and passing if necessary to a subsequence
$$w_{a}\rightharpoonup w_0\,\ \hbox{weakly in} \,\ H^1(\R^N),\,\  w_{a}\rightarrow w_0\,\ \hbox{strongly in}\,\ L_{loc}^{q}(\R^{N}),\,\ 2\leq q<2^*. $$
Note that the proof of the claim $\eqref{2.2}$ do not rely on the assumption $V(0)>0$ and $e(a)< V(0)$. Thus, the same argument of proving \eqref{2.2} yields that
\begin{equation}\label{2.5}
w_{a}\rightarrow \frac{\beta^{\frac{N}{2}}Q(\beta x)}{\|Q\|_2}\,\ \hbox{strongly in}\,\ H^1(\R^N)\,\ \hbox{as}\,\ a\nearrow a^*.
\end{equation}

Next we prove that
\begin{equation}\label{2.6}
w_a(x)\rightarrow \frac{\beta^{\frac{N}{2}}Q(\beta x)}{\|Q\|_2}\,\ \hbox{strongly in}\,\ L^\infty(\R^N)\,\ \hbox{as}\,\  a\nearrow a^*.
\end{equation}
For $N=1$, this is trivial since the embedding $H^1(\R)\hookrightarrow L^\infty(\R)$ is continuous.

We  next focus on the case $N\geq2$.
We claim that
\begin{equation}\label{21}
w_{a}(x)\rightarrow 0 \ \ \mbox{as}\ \ |x|\rightarrow\infty\ \ \mbox{uniformly}\ \ \mbox{for}\ \  a\nearrow a^*.
\end{equation}
In fact, we derive from \eqref{2.6} that for $2\leq q< 2^*$
\begin{equation}\label{22}
\int_{|x|\geq \gamma}|w_{a}(x)|^q\rightarrow 0 \ \ \mbox{as}\ \ \gamma\rightarrow\infty\ \ \mbox{uniformly}\ \ \mbox{for}\ \  a\nearrow a^*.
\end{equation}
On the other hand, recall from \eqref{EL} and \eqref{xpa} that $w_a$ satisfies the following equation
\begin{equation}\label{2.9}
-\Delta w_a+\eps_a^2V(\eps_ax)w_a=\mu_a\eps_a^2w_a+a w_a^{1+2\beta^2} |x|^{-b}\quad \hbox{in}\,\ \R^N,
\end{equation}
where $\mu_a\in\R$ is a suitable Lagrange multiplier  satisfying \eqref{2.8}.
In view of \eqref{2.4} and \eqref{2.8}, we derive from the convergence in \eqref{2.6} that
\begin{equation}\label{2.10}
\begin{split}
\eps_a^2\mu_a&=\eps_a^2e(a)-\eps_a^2\frac{a\beta^2}{1+\beta^2}\inte |x|^{-b}|u_a|^{2+2\beta^2}dx
\\
&=\eps_a^2e(a)-\frac{a\beta^2}{1+\beta^2}\inte |x|^{-b}|w_a|^{2+2\beta^2}dx\to -\beta^2\,\ \hbox{as}\,\ a\nearrow a^*.
\end{split}
\end{equation}
Hence for $a\nearrow a^*$, we have
\begin{equation}\label{R1}
-\Delta w_{a}-c(x)w_{a}(x)\leq0\ \ \mbox{in}\ \ \R^N,\ \ \mbox{where}\ \ c(x)=w^{2\beta^2}_{a}|x|^{-b}.
\end{equation}
Note from Lemma \ref{A.1} in the Appendix that $c(x)\in L^{m}(\R^N)$ where $m\in \Big(\frac{N}{2},\frac{N^2}{2N+2b-4}\Big)$ if $N\geq 3$ and $m\in(1,\frac 2b)$ if $N=2$.
Applying De-Giorgi-Nash-Moser theory  \cite[Theorem 4.1]{QF} to \eqref{R1}, we conclude that
\begin{equation}\label{221}
\max\limits_{B_1(\rho)}w_{a}(x)\leq C\Big(\int_{{B_2(\rho)}}|w_{a}(x)|^qdx\Big)^{\frac{1}{q}}\ \ \mbox{as}\ \ a\nearrow a^*.
\end{equation}
where $\rho\in\R^N$ is arbitrary, and $C>0$ depends only on the bound of $\|c(x)\|_{L^{m}(\R^N)}$. Therefore, the claim \eqref{21} is  proved in view of \eqref{22} and \eqref{221}.

Based on \eqref{21}, to finish the proof of \eqref{2.6}, the rest is to show that
\begin{equation}\label{2.11}
w_a\to \frac{\beta^{\frac{N}{2}}Q(\beta x)}{\|Q\|_2}\,\ \hbox{strongly in}\,\ L^\infty_{loc}(\R^N).
\end{equation}
Since $w_a(x)$ satisfies (\ref{2.9}), we denote
\begin{equation*}
G_a(x):=\varepsilon_a\mu_aw_a(x)-\varepsilon_a^2V(\varepsilon_ax)w_a(x)+w^{1+2\beta^2}_{a}|x|^{-b},
\end{equation*}
so that
\begin{equation}\label{g1}
-\Delta w_a(x)=G_a(x) \ \ \mbox{in}\ \ H^1(\R^N).
\end{equation}
Because $w_a(x)$ is bounded uniformly in $H^1(\R^N)$, follow (\ref{221}), we can obtain that
\begin{equation}\label{ubl}
w_a\ \ \mbox{is bounded}\ \ \mbox{uniformly}\ \ \mbox{in} \ \ L^\infty(\R^N).
\end{equation}
Therefore, for any $r\in(1,\frac{N}{b})$, $w_{a}^{1+2\beta^2}|x|^{-b}$ is bounded uniformly in $L^r_{loc}(\R^N)$, which implies that $G_a(x)$ is bounded uniformly in $L^r_{loc}(\R^N)$. For any large $R>0$, it then follows from \cite[Theorem 9.11]{GT} that
\begin{equation}\label{wgj}
\|w_a(x)\|_{W^{2,r}(B_R)}\leq C\Big(\|w_a(x)\|_{L^r(B_{R+1})}+\|G_a(x)\|_{L^r(B_{R+1})}\Big),
\end{equation}
where $C>0$ is independent of $a>0$ and $R>0$. Therefore, $w_a(x)$ is also uniformly in $W^{2,r}_{loc}(\R^N)$. Since for $r>\frac{N}{2}$, the embedding $W^{2,r}(B_R)\hookrightarrow L^\infty(B_R)$ is compact, see \cite[Theorem 7.26]{GT}, we then deduce that there exists a subsequence $\{w_{a_{k}}\}$ of $\{w_a(x)\}$ such that
\begin{equation*}
\lim\limits_{a_k\nearrow a^*}w_{a_k}(x)=\tilde{w}_0(x)\ \ \mbox{uniformly}\ \ \mbox{in}\ \ L^\infty(B_R).
\end{equation*}
In view of (\ref{2.5}) and the fact that $R>0$ is arbitrary, we obtain that
\begin{equation}\label{lwj2}
\lim\limits_{a_k\nearrow a^*}w_{a_k}(x)=\frac{\beta^{\frac{N}{2}}Q(\beta x)}{\|Q\|_2}\ \ \mbox{uniformly}\ \ \mbox{in}\ \ L_{loc}^\infty(\R^N).
\end{equation}
Since the above convergence is independent of what subsequence we choose,  we deduce that $w_{a}(x)\rightarrow\frac{\beta^{\frac{N}{2}}Q(\beta x)}{\|Q\|_2}$ in $L_{loc}^{\infty}(\R^N)$ as $a \nearrow a^*$ and (\ref{2.6}) is proved in view of \eqref{21}.

(3).  Following \eqref{21}, \eqref{2.9} and \eqref{2.10}, we obtain that there exists $R>0$ large enough such that as $a\nearrow a^*$
\begin{equation}\label{bj}
-\Delta w_{a}(x)+\frac{\beta^2}{4}w_{a}(x)\leq 0\ \ \mbox{in}\ \ \R^N\backslash B_R(0).
\end{equation}
Then by the comparison principle, we then have as $a \nearrow a^*$,
\begin{equation}\label{bj1}
w_{a}(x)\leq Ce^{-\frac{\beta|x|}{2}}\ \ \mbox{in}\ \ \R^N \backslash B_R(0).
\end{equation}
Combing the above estimate and the fact (\ref{ubl}), we conclude that
\begin{equation}\label{qbj1}
w_{a}(x)\leq Ce^{-\frac{\beta|x|}{2}}\ \ \mbox{in}\ \ \R^N.
\end{equation}
Moreover, applying the local elliptic estimate \cite[(3.15)]{GT} to \eqref{g1} yields that as $a\nearrow a^*$,
\begin{equation}\label{bj2}
|\nabla w_{a}(x)|\leq Ce^{-\frac{\beta|x|}{4}}\ \ \mbox{in}\ \ \R^N.
\end{equation}
Therefore, the exponential decay (\ref{zsll}) is then proved.  \qed

\vskip 0.1truein
\noindent\textbf{Completion the proof of Theorem 1.3.} In view of Lemma \ref{yl1}, to establish Theorem \ref{th3}, it remains to prove that
\begin{equation}\label{xpx111}
\begin{split}
\varepsilon_{a}&:=\Big(\int_{\R^N}|\nabla u_a|^2\Big)^{-\frac{1}{2}}
\\&=\frac{(a^{\ast}-a)^{\frac{1}{l+2}}}{\lambda}+o\big((a^*-a)^{\frac{1}{l+2}}\big)\ \ \mbox{as}\ \ a\nearrow a^*,
\end{split}
\end{equation}
where $\lambda>0$ is as in Theorem \ref{th3}.

In order to prove (\ref{xpx111}), we take the following test function
\begin{equation*}
u_{\alpha}(x):=\frac{\alpha^{\frac{N}{2}}}{\|Q\|_{2}}Q (\alpha x),
\end{equation*}
where $\alpha\in (0, \infty)$ is to be determined later. Direct calculations give that as $\alpha\to\infty$
\begin{equation*}
\int_{\R^{N}}|\nabla u_{\alpha}(x)|^{2}dx=\frac{\alpha^{2}}{\|Q\|^{2}_{2}}\int_{\R^{N}}|\nabla Q|^{2}dx=\frac{\alpha^2}{\beta^2},
\end{equation*}
\begin{equation*}
\begin{split}
&\frac{a}{1+\beta^2}\int_{\R^{N}}|x|^{-b}|u_{\alpha}(x)|^{2+2\beta^2}dx\\
=&\frac{a\alpha^{2}}{(1+\beta^2)\|Q\|_2^{2+2\beta^2}}\int_{\R^{N}}|x|^{-b}|Q(x)|^{2+2\beta^2}dx=\frac{a}{a^*}\frac{\alpha^2}{\beta^2},
\end{split}
\end{equation*}
and
\begin{equation*}
\int_{\R^{N}}V(x)u^{2}_{\alpha}(x)dx=\frac{1}{\|Q\|^{2}_{2}}\int_{\R^{N}}V\Big(\frac{x}{\alpha}\Big)Q^{2}(x)dx= \frac{1}{\alpha^{l}\|Q\|^{2}_{2}}\int_{\R^{N}}h(x)Q^{2}(x)dx+o(\alpha^{-l}).
\end{equation*}
Thus we have
\begin{equation}\label{nlxj1}
\begin{split}
e(a)\leq E_{a}(u_{\alpha})&=\frac{\alpha^{2}}{\beta^{2}}\frac{a^{\ast}-a}{a^{\ast}}+\frac{1}{\alpha^{l}\|Q\|^{2}_{2}}\int_{\R^{N}}h(x)Q^{2}(x)dx+o(\alpha^{-l})
\\&:= g(\alpha)+o(\alpha^{-l}).
\end{split}
\end{equation}
Therefore, we have
\begin{equation}\label{xea}
\begin{split}
e(a)&\leq \min_{\alpha\in(0,+\infty)}g(\alpha)+o(\alpha^{-l})=g\Big(\Big(\frac{\int_{\R^{N}}h(x)Q^{2}(x)dx\beta^2la^*}{2\|Q\|_2^2(a^*-a)}\Big)^{\frac{1}{l+2}}\Big)+o((a^*-a)^{\frac{l}{l+2}})\\
&=\Big( \frac{\int_{\R^{N}}h(x)Q^{2}(x)dx}{\|Q\|^{2}_{2}\beta^{l}}\Big)^{\frac{2}{l+2}}\Big[\Big(\frac{l}{2}\Big)^{\frac{2}{l+2}}+\Big(\frac{2}{l}\Big)^{\frac{l}{l+2}}\Big]\Big(\frac{a^{\ast}-a}{a^{\ast}}\Big)^{\frac{l}{l+2}}+o((a^*-a)^{\frac{l}{l+2}}).
\end{split}
\end{equation}
On the other hand,
\begin{equation}\label{xnlgj}
\begin{split}
e(a)&=E_{a}(u_{a})
\\&=\int_{\R^{N}}|\nabla u_{a}(x)|^{2}dx-\frac{a}{1+\beta^2}\int_{\R^{N}}|x|^{-b}u_{a}^{2+2\beta^2}(x)dx+\int_{\R^{N}}V(x)|u_{a}(x)|^{2}dx
\\&=\frac{1}{\varepsilon^{2}_{a}}\Big(\int_{\R^{N}}|\nabla w_{a}(x)|^{2}dx-\frac{a^{\ast}}{1+\beta^2}\int_{\R^{N}}|x|^{-b}w_{a}^{2+2\beta^2}(x)dx\Big)
\\&+\frac{a^{\ast}-a}{(1+\beta^2)\varepsilon^{2}_{a}}\int_{\R^{N}}|x|^{-b}w_{a}^{2+2\beta^2}(x)dx+\int_{\R^{N}}V(\varepsilon_{a}x)|w_{a}(x)|^{2}dx.
\end{split}
\end{equation}
Under the assumption $(V_2)$, we deduce from Lemma \ref{yl1} (2) and (3) that
\begin{equation}\label{va}
\begin{split}
&\int_{\R^{N}}V(\varepsilon_{a}x)|w_{a}(x)|^{2}dx
\\
&=\int_{B_{\varepsilon_{a}^{-\frac 12}}(0)}V(\varepsilon_{a}x)|w_{a}(x)|^{2}dx+\int_{\R^{N}/{B_{\varepsilon_{a}^{-\frac 12}}(0)}}V(\varepsilon_{a}x)|w_{a}(x)|^{2}dx
\\
&=\varepsilon^{l}_{a}\int_{B_{\varepsilon_{a}^{-\frac 12}}(0)}h(x)|w_{a}(x)|^{2}dx+o(\varepsilon^{l}_{a})
\\
&=\frac{[1+o(1)]\varepsilon^{l}_{a}}{\|Q\|_2^{2}\beta^l}\int_{\R^{N}}h(x)|Q(x)|^{2}dx+o(\varepsilon^{l}_{a})\ \ \mbox{as}\ \ a \nearrow a^*.
\end{split}
\end{equation}
Similarly, we deduce from \eqref{2.10} that
\begin{equation}\label{slsl}
\begin{split}
&\frac{a^{\ast}-a}{(1+\beta^2)\varepsilon^{2}_{a}}\int_{\R^{N}}|x|^{-b}w_{a}^{2+2\beta^2}(x)dx=[1+o(1)]\frac{(a^*-a)}{a^*\varepsilon^{2}_{a}}\ \ \mbox{as}\ \ a \nearrow a^*.
\end{split}
\end{equation}
Therefore,  we obtain that as $a\nearrow a^*$
\begin{equation}\label{tnli}
\begin{split}
e(a)&\geq \frac{a^{\ast}-a}{a^{\ast}\varepsilon^{2}_{a}}+\frac{\varepsilon^{l}_{a}}{\beta^{l}\|Q\|^{2}_{2}}\int_{\R^{N}}h(x)|Q(x)|^{2}dx+o(\varepsilon^{l}_{a}+\frac{a^{\ast}-a}{\varepsilon^{2}_{a}})\\
&=g(\beta\varepsilon_a^{-1})+o(\varepsilon^{l}_{a}+\frac{a^{\ast}-a}{\varepsilon^{2}_{a}})\\
&\geq \min_{t\in(0,+\infty)}g(t)+o((a^*-a)^{\frac{l}{l+2}})\\
&=\Big( \frac{\int_{\R^{N}}h(x)Q^{2}(x)dx}{\|Q\|^{2}_{2}\beta^{l}}\Big)^{\frac{2}{l+2}}\Big[\Big(\frac{l}{2}\Big)^{\frac{2}{l+2}}+\Big(\frac{2}{l}\Big)^{\frac{l}{l+2}}\Big]\Big(\frac{a^{\ast}-a}{a^{\ast}}\Big)^{\frac{l}{l+2}}+o((a^*-a)^{\frac{l}{l+2}}).
\end{split}
\end{equation}
Note from the upper bound estimates of energy in \eqref{xea} that the above inequalities must be achieved and hence $g(\beta\eps_a^{-1})=\min_{t\in(0,+\infty)}g(t)$. This further implies that
\begin{equation*}
\varepsilon_{a}=\frac{(a^{\ast}-a)}{\lambda}^{\frac{1}{l+2}}+o\big((a^*-a)^{\frac{1}{l+2}}\big)\quad\hbox{as}\,\ a\nearrow a^*,
\end{equation*}
and the proof of Theorem \ref{th3} is thus completed.\qed

\section{Proof of local uniqueness}
This section is devoted to  the proof of Theorem \ref{unique} on the uniqueness of  positive minimizers as $a\nearrow a^*$. Stimulated by \cite{Cao,Deng},
we first define the linearized operator $\mathcal{L}$ by
\[
\mathcal{L}:=-\Delta +1-(1+2\beta^2)|x|^{-b}Q^{2\beta^2}\ \  \mbox{in}\ \, \R^N,
\]
where $Q=Q(|x|)>0$ is the unique positive  solution of (\ref{Eq}). Note from \cite{SC} that
\begin{equation}
ker (\mathcal{L})=span \{0\}\quad\hbox{if}\,\ N\geq 3.
\label{2:linearized}
\end{equation}
Recall that $u_a$ solves the Euler-Lagrange equation
\begin{equation}
-\Delta u_a(x) +V(x)u_a(x) =\mu_a u_a(x) +a u ^{1+2\beta^2} _a(x)|x|^{-b}\ \, \text{in\  $\R^N$},
\label{2:cond}
\end{equation}
where $\mu _a\in \R$ is a suitable Lagrange multiplier and satisfies
\begin{equation}\label{2:cone}
\mu_a = e(a) - \frac{\beta^2a}{1+\beta^2} \int_{\R^{N}} u^{2(1+\beta^2)}_a|x|^{-b} dx.
\end{equation}
For simplicity, in the whole section we always denote
\begin{equation}\label{4.1}
\alp_a:=\frac{(a^*-a)^{ \frac{1}{2+l}}}{\lambda}>0.
\end{equation}
Following \eqref{xpx111}, we have
\begin{equation}\label{5.4}
\lim_{a\nearrow a^*}\frac{\eps_a}{\alp_a}=1,
\end{equation}
which combined with \eqref{2.10} gives that
\begin{equation}\label{uyp}
\mu_{a}\alp_{a}^{2}\rightarrow -\beta^{2}\,\ \mbox{as} \,\ a\nearrow a^*.
\end{equation}

\vskip 0.1truein
\noindent\textbf{Proof of Theorem \ref{unique}.}
Suppose that there exist two different positive minimizers $u_{1,a}$ and $u_{2,a}$ of $e(a)$ as $a\nearrow a^*$.  Following (\ref{2:cond}), $u_{i,a}$ then solves the Euler-Lagrange equation
\begin{equation}
-\Delta u_{i,a}(x) +V(x)u_{i,a}(x) =\mu_{i,a} u_{i,a}(x) +au^{1+2\beta^2} _{i,a}(x)|x|^{-b}\quad \text{in\ \ $\R^N$}, \ \ i=1,2,
\label{3d:cond}
\end{equation}
where $\mu _{i,a}\in \R$ is a suitable Lagrange multiplier.
Define
\begin{equation}\label{3d:27A}
\bar u_{i,a}(x):=\frac{\alp_{a}^{\frac{N}{2}}\|Q\|_{2}}{\beta^{\frac{N}{2}}}u_{i,a} ( x ), \ \ \mbox{where}\ \ i=1,2.
\end{equation}
Then Theorem \ref{th3}, Lemma \ref{yl1} (3) and \eqref{5.4} imply that for $i=1,2,$
\begin{equation}\label{5.5}
\bar u_{i,a}\big(\frac{\alp_a}{\beta}x\big)\to Q(x)\quad\hbox{stronly in $L^{\infty}(\R^{N})\cap H^1(\R^N)$ as $a\nearrow a^*$.}
\end{equation}
and
\begin{equation}\label{na:A}
\big|\bar u_{i,a}\big(\frac{\alp_a}{\beta}x\big)\big|,\, \ \big|\nabla \bar u_{i,a}\big(\frac{\alp_a}{\beta}x\big)\big|\le C e^{-\frac{|x|}{4}}\,\ \mbox{in}  \,\ \R^N,
\end{equation}
where $C>0$ is independent of $a\nearrow a^*$.
For $i=1,2$, $\bar u_{i,a}$ also satisfies the equation
\begin{equation}
-\alp ^2 _a\Delta \bar u_{i,a}(x) +\alp ^2 _aV(x)\bar u_{i,a}(x) =\mu_{i,a}\alp ^2 _a\bar  u_{i,a}(x) +\frac{\alp^{b}_{a}\beta^{2-b}a}{a^*}\bar u^{1+2\beta^2} _{i,a}(x)|x|^{-b}\quad \text{in\ \ $\R^N$}.
\label{5-2:1}
\end{equation}
Because $u_{1,a}\not\equiv u_{2,a}$, we consider
\begin{equation}\label{5.6}
\bar \xi_a(x)=\frac{ u_{2,a}(x)-  u_{1,a}(x)}{\|  u_{2,a}-  u_{1,a}\|_{L^\infty(\R^N)}}=\frac{\bar u_{2,a}(x)-\bar u_{1,a}(x)}{\|\bar u_{2,a}-\bar u_{1,a}\|_{L^\infty(\R^N)}}.
\end{equation}
Motivated by \cite{Cao}, we first claim that for any $x_0\in\R^N$, there exists a small constant $\delta >0$  such that
\begin{equation}
\int_{\partial B_\delta (x_0)} \Big[ \alp ^2_a |\nabla \bar \xi_a|^2+ \frac{\beta ^2}{2} |\bar\xi_a|^2+ \alp ^2_a  V(x)|\bar\xi_a|^2\Big]dS=O( \alp ^N_a)\quad \text{as}\ a\nearrow a^*.
\label{5.2:6}
\end{equation}
Since the proof of the claim is standard and lengthy (see \cite{Cao,GLW}), we left its to the Appendix.

We next define
\begin{equation}\label{5.1}
\tilde{u}_{i,a}(x):=\bar{u}_{i,a} \big(\frac{\alp _a}{\beta}x\big), \ \ \mbox{where}\ \ i=1,2,\ \ \hbox{as}\,\ a\nearrow a^*,
\end{equation}
and
\begin{equation}
\xi _a(x)=\bar \xi_a\big(\frac{\alp _a}{\beta}x\big) \ \ \hbox{as}\,\ a\nearrow a^*,
\label{5.2:xi}
\end{equation}
then   $\tilde{u}_{i,a}(x)\to Q(x)$ uniformly in $\R^N$ as $a\nearrow a^*$ in view of \eqref{5.5}. We shall carry out the proof of Theorem \ref{unique} by deriving a contradiction through the following three steps.

\vskip 0.1truein

\noindent{\em  Step 1.} There exist a subsequence of $\{a\}$, still denote by $\{a\}$,  and constant $b_0\in\R$ such that $\xi_a(x)\to \xi_0(x)$ in $C_{loc}(\R^N)$ as $a\nearrow a^*$, where
\begin{equation}\label{ks0}
\xi_{0}=b_0\big(\frac{N}{2}Q+x\cdot\nabla Q\big).
\end{equation}

Indeed, note that  $\tilde{u}_{i,a}(x)$ satisfies the following equation:
\begin{equation}\label{xy}
-\Delta \tilde{u}_{i,a}(x)+\frac{\alp^{2}_{a}}{\beta^{2}}V\big(\frac{\alp_{a}}{\beta}x\big)\tilde{u}_{i,a}(x)=\frac{a}{a^{\ast}}|x|^{-b} \tilde{u}_{i,a}^{1+2\beta^2}(x)+ \frac{\alp^{2}_{a}}{\beta^{2}}\mu_{i,a}\tilde{u}_{i,a}(x).
\end{equation}
We denote
\begin{equation}\label{dkp2}
\arraycolsep=1.5pt\begin{array}{lll}
\tilde{D}^{s-1}_{a}(x):=\int^{1}_{0}\big[t\tilde{u}_{2,a}+(1-t)\tilde{u}_{1,a}\big]^{s-1}dt,
\end{array}
\end{equation}
then $\xi_a$ satisfies
\begin{equation}
-\Delta \xi_a +C_{a}(x)\xi_a =g_a(x)\quad \text{in\,\, $\R^N$},
\label{3d:2}
\end{equation}
where the coefficient $C_a(x)$ satisfies
\begin{equation}\label{4.20}
C_a(x):=-\frac{a(1+2\beta^2)\tilde{D}^{2\beta^2}_{a}(x)}{a^*|x|^b}-\frac{\alp ^2_a}{\beta ^2}\mu _{1,a}+\frac{\alp ^2_a}{\beta ^2}V\big(\frac{\alp_ax}{\beta}\big),
\end{equation}
and the nonhomogeneous term $g_a(x)$ satisfies
\begin{equation}
\arraycolsep=1.5pt\begin{array}{lll}
g_a(x)&:=&\displaystyle\frac{\tilde{u}_{2,a}}{\beta ^2}\frac{\alp ^2_a(\mu _{2,a}-\mu _{1,a})}{\|\tilde{u}_{2,a}-\tilde{u}_{1,a}\|_{L^\infty}}\\[4mm]
&=&-\displaystyle \frac{2a\beta^2\tilde{u}_{2,a}}{\|Q\|_2^{2(1+\beta^2)}}  \int_{\R^{N}} |x|^{-b}\xi_{a}\tilde{D}^{1+2\beta^2}_{a}dx.
\end{array}
\label{3d:4}
\end{equation}
Here we have used (\ref{2:cone}) and (\ref{5.1}). Recall that $\tilde u_{i,a}$ is bounded uniformly in $L^\infty(\R^N)$, we obtain that $|x|^{-b}\tilde{D}_a^{2\beta^2}$ is bounded uniformly $L^r_{loc}(\R^N)$, where $r\in(1,\frac{N}{b})$. 
Since $\|\xi _a\|_{L^\infty(\R^N)}\le 1$, the standard elliptic regularity then implies cf.(\cite{GT}) that $\|\xi _a\|_{W^{2,r}_{loc}(\R^N)}\le C$  and thus $\|\xi _a\|_{C^{\alpha }_{loc}(\R^N)}\le C$ for some $\alp \in (0,2-b)$, where the constant $C>0$ is independent of $a$. Therefore, there exists a subsequence of $\{a\}$ still denote by $\{a\}$ and a function $\xi _0=\xi _0(x)$ such that $\xi_a(x)\to \xi_0(x)$ in $C_{loc}(\R^N)$ as $a\nearrow a^*$.
Applying Lemma \ref{yl1}, direct calculations yield from (\ref{4.20}) and \eqref{3d:4}  that
\[
C_a(x)\to 1-(1+2\beta^2)|x|^{-b}Q^{2\beta^2}(x)\ \ \text{uniformly\ in\ $C_{loc}(\R^N)$}\,\ \mbox{as} \,\ a\nearrow a^*,
\]
and
\[
g_a(x)\to -\frac{2\beta^2Q(x)}{\|Q\|^{2}_{2}}\int_{\R^{N}}|x|^{-b}\xi _0Q^{1+2\beta^2}\ \ \text{uniformly\ in\ $C_{loc}(\R^N)$}\,\ \mbox{as} \,\ a\nearrow a^*.
\]
This implies from (\ref{3d:2}) that $\xi_0$ solves
\begin{equation}
\mathcal{L}\xi_0=-\Delta \xi_0+\xi_0-(1+2\beta^2)|x|^{-b}Q^{2\beta^2}\xi_0=-\frac{2\beta^2Q(x)}{\|Q\|^{2}_{2}}\int_{\R^{N}}|x|^{-b}\xi _0Q^{1+2\beta^2}\ \ \mbox{in} \ \ \R^N.
\label{3d:5}
\end{equation}
Since $\mathcal{L}\big(\frac{N}{2}Q+x\nabla Q\big)=-2Q$, see Lemma \ref{A.2} in the Appendix, we then derive from the non-degeneracy of $\mathcal{L}$ in \eqref{2:linearized} that (\ref{ks0}) holds true for a certain constant $b_0$.
\vskip 0.1truein

\noindent{\em  Step 2.} The constant $b_0=0$ .
Using the integration by parts,  we obtain that
\begin{equation}\arraycolsep=1.5pt\begin{array}{lll}
&&-\displaystyle\alp _a^2 \int_{B_{\delta}(0)}(x\cdot \nabla \bar u_{i,a}) \Delta \bar u_{i,a} \\[4mm]
&=&-\alp _a^2\displaystyle  \int_{\partial B_{\delta}(0)} \frac{\partial\bar u_{i,a} }{\partial \nu }(x\cdot \nabla \bar u_{i,a})dS
+\displaystyle\alp _a^2 \int_{B_{\delta}(0)} \nabla \bar u_{i,a}\nabla (x\cdot \nabla \bar u_{i,a}) \\[4mm]
&=&\alp _a^2\displaystyle  \int_{\partial B_{\delta}(0)} \Big[-\frac{\partial\bar u_{i,a} }{\partial \nu }(x\cdot \nabla \bar u_{i,a})+\frac 12(x\cdot \nu) |\nabla \bar u_{i,a}|^2\Big]dS
+\frac{(2-N)\alp_a^2}{2}\int_{B_{\delta}(0)}|\nabla \bar u_{i,a}|^2\\[4mm]
&=&\alp _a^2\displaystyle  \int_{\partial B_{\delta}(0)} \Big[-\frac{\partial\bar u_{i,a} }{\partial \nu }(x\cdot \nabla \bar u_{i,a})+\frac 12(x\cdot \nu) |\nabla \bar u_{i,a}|^2+\frac{(2-N)}{4}(\nabla \bar u_{i,a}^2\cdot\nu) \Big]dS\\[4mm]
&&\quad-\displaystyle\frac{(2-N)}{2} \int_{B_{\delta}(0)}\Big[\alp _a^2V(x)\bar u^{2}_{i,a}-\alp _a^2\mu_{i,a}\bar u_{i,a}^{2}- \frac{\alp _a^b\beta^{2-b}a\bar u_{i,a}^{2(1+\beta^2)}}{a^*|x|^b}\Big],
\end{array}\label{5.3:3}
\end{equation}
where the last equality used the following fact
\begin{equation}\arraycolsep=1.5pt\begin{array}{lll}
&&\displaystyle\frac{(2-N)\alp _a^2}{2}\int_{B_{\delta}(0)} |\nabla \bar u_{i,a}|^2
=\displaystyle\frac{(2-N)\alp _a^2}{4}\int_{\partial B_{\delta}(0)}(\nabla \bar u_{i,a}^2\cdot\nu) dS\\
&&\displaystyle\qquad-\frac{2-N}{2} \int_{B_{\delta}(0)}\Big[\alp _a^2V(x)\bar u^{2}_{i,a}-\alp _a^2\mu_{i,a}\bar u_{i,a}^{2}- \frac{\alp _a^b\beta^{2-b}a\bar u_{i,a}^{2(1+\beta^2)}}{a^*|x|^b}\Big].
\end{array}\label{5.4:2}
\end{equation}

On the other hand, multiplying (\ref{5-2:1}) by $ (x\cdot \nabla  \bar{u}_{i,a})$, where $i=1,2$, and integrating over $B_\delta (0)$, where $\delta >0$ is small as before, we deduce that for $i=1,2,$
\begin{equation}\arraycolsep=1.5pt\begin{array}{lll}
&&-\displaystyle\alp _a^2 \int_{B_{\delta}(0)} \big(x\cdot \nabla \bar u_{i,a}\big) \Delta \bar u_{i,a} \\[4mm]
&=& \displaystyle\alp _a^2 \int_{B_{\delta}(0)} \big[\mu_{i,a}-V(x)\big] \bar u_{i,a}\big(x\cdot \nabla \bar u_{i,a}\big)+
\displaystyle \frac{\beta^{2-b}\alp^{b}_{a} a}{a^*} \int_{B_{\delta}(0)} \bar u_{i,a}^{1+2\beta^2}|x|^{-b}\big (x\cdot \nabla \bar u_{i,a}\big)\\[4mm]
&=& -\displaystyle\frac{\alp _a^2}{2} \int_{B_{\delta}(0)}\bar u_{i,a}^2\Big\{N\big[\mu_{i,a}-V(x)\big]-[x\cdot \nabla V(x)]\Big\}+\displaystyle\frac{\alp _a^2}{2}  \int_{\partial B_{\delta}(0)}\bar u_{i,a}^2\big[\mu_{i,a}-V(x)\big](x\cdot\nu) dS\\[4mm]
&&+\displaystyle \frac{\beta^{2-b}\alp^{b}_{a} a}{2(1+\beta^2)a^*} \Big[\displaystyle \int_{\partial B_{\delta}(0)} \bar{u}_{i,a}^{2(1+\beta^2)}\frac{(x\cdot\nu)}{|x|^{b}} dS-\int_{B_{\delta}(0)} \bar u_{i,a}^{2(1+\beta^2)}\frac{N-b}{|x|^{b}}\Big]\\[4mm]
\end{array}\label{5.3:1}
\end{equation}
Substituting (\ref{5.3:3}) into (\ref{5.3:1}) yields that
\begin{equation}\arraycolsep=1.5pt\begin{array}{lll}
&&\alp _a^2\displaystyle  \int_{\partial B_{\delta}(0)} \Big[-\frac{\partial\bar u_{i,a} }{\partial \nu }(x\cdot \nabla \bar u_{i,a})+\frac 12(x\cdot \nu) |\nabla \bar u_{i,a}|^2+\frac{(2-N)}{4}(\nabla \bar u_{i,a}^2\cdot\nu) \Big]dS\\[4mm]
&=&\displaystyle \int_{B_{\delta}(0)}\Big\{\alp _a^2\big[-\mu_{i,a}+V(x)+\frac{1}{2}(x\cdot\nabla V)\big] \bar u_{i,a}^2-\frac{\beta^{4-b}\alp^{b}_{a} a}{(1+\beta^2)a^*}\bar u_{i,a}^{2(1+\beta^2)}|x|^{-b}\Big\}+I_i\\[4mm]
\end{array}\label{5.4:3}
\end{equation}
where the lower order term $I_i$ satisfies
\begin{equation}\arraycolsep=1.5pt\begin{array}{lll}
I_i&=&\displaystyle\frac{\alp _a^2}{2} \int_{\partial B_{\delta}(0)}\bar u_{i,a}^2\big[\mu_{i,a}-V(x)\big](x\cdot\nu) dS+\displaystyle \frac{\beta^{2-b}\alp^{b}_{a} a}{2(1+\beta^2)a^*} \int_{\partial B_{\delta}(0)} \bar u_{i,a}^{2(1+\beta^2)}\frac{(x\cdot\nu)}{|x|^{b}}dS
,\,\ i=1,2.
\end{array}\label{5.3:2}
\end{equation}
Since it follows from (\ref{2:cone}) that
\[\arraycolsep=1.5pt\begin{array}{lll}
&&\quad\displaystyle \mu_{i,a}\displaystyle\alp _a^2 \int_{\R^{N}} \bar u_{i,a}^2+\displaystyle \frac{\beta^{4-b}\alp^{b}_{a} a}{(1+\beta^2)a^*} \int_{\R^{N}}  \frac{\bar u_{i,a}^{2(1+\beta^2)}}{|x|^{b}}=\displaystyle\frac{\alp^{N+2}_{a}\|Q\|^{2}_{2}}{\beta^{N}}e(a),
\end{array}\]
we reduce from (\ref{5.3:3})--(\ref{5.3:2}) that
\[\arraycolsep=1.5pt\begin{array}{lll}
&&\quad\displaystyle-\alp^{2}_{a}\int_{B_{\delta}(0)} \Big[V(x)+\frac{1}{2}[x\cdot \nabla V(x)]\Big]\bar u_{i,a}^{2}+\frac{\alp^{N+2}_{a}\|Q\|^{2}_{2}}{\beta^{N}}e(a)\\[4mm]
&&=I_i+\alp _a^2\displaystyle \int_{\partial B_{\delta}(0)} \frac{\partial\bar u_{i,a} }{\partial \nu }(x\cdot \nabla \bar u_{i,a})
-\displaystyle \frac{\alp _a^2}{2}\int_{\partial B_{\delta}(0)} (x\cdot \nu )|\nabla \bar u_{i,a}|^2\\[4mm]
&&\quad-\displaystyle\frac{2-N}{4}\alp _a^2\int_{\partial B_{\delta}(0)} \nabla {\bar u_{i,a}}^{2}+
\mu_{i,a}\alp _a^2 \int _{\R^N\backslash B_\delta (0)} \bar u_{i,a}^2\\[4mm]
&&\quad+\displaystyle \frac{\beta^{4-b}\alp^{b}_{a} a}{(1+\beta^2)a^*} \int _{\R^N\backslash B_\delta (0)} \bar u_{i,a}^{2(1+\beta^2)}|x|^{-b},
\,\ i=1,2,
\end{array}\]
which implies that
\begin{equation}
\displaystyle-\alp^{2}_{a}\int_{B_{\delta}(0)} \Big[V(x)+\frac{1}{2}[x\cdot \nabla V(x)]\Big](\bar u_{1,a}+\bar u_{2,a})\bar{\xi}_a=T_a.
\label{5.3:4}
\end{equation}
Here the term $T_a$ satisfies that for small $\delta >0$,
\begin{equation}\arraycolsep=1.5pt\begin{array}{lll}
T_a&=&\displaystyle\frac{I_2-I_1}{\|\bar u_{2,a}-\bar u_{1,a}\|_{L^\infty}}-\displaystyle \frac{2-N}{4}\alp _a^2\int_{\partial B_{\delta}(0)}\big(\nabla \bar u_{2,a}+\nabla \bar u_{1,a}\big)\nabla \bar \xi_a\\[4mm]
&&-\displaystyle\frac{\alp _a^2}{2}\int_{\partial B_{\delta}(0)} (x\cdot \nu )\big(\nabla \bar u_{2,a}+\nabla \bar u_{1,a}\big)\nabla \bar \xi_a\\[4mm]
&&+\displaystyle  \alp _a^2 \int_{\partial B_{\delta}(0)} \Big[(x\cdot \nabla \bar u_{2,a})\big(\nu \cdot \nabla\bar \xi_a \big)+\big(\nu \cdot \nabla \bar u_{1,a}\big)(x \cdot \nabla \bar \xi_a)\Big]
\\[4mm]
&&+\mu_{2,a}\displaystyle\alp _a^2 \int _{\R^N\backslash B_\delta (0)} \big(\bar u_{2,a}+\bar u_{1,a}\big)\bar \xi_a+\frac{(\mu_{2,a}-\mu_{1,a})\alp_a^2}{\|\bar u_{2,a}-\bar u_{1,a}\|_{L^\infty}}\int _{\R^N\backslash B_\delta (0)} \bar u_{1,a}^2\\[4mm]
&&+\displaystyle\frac{\alp^{b}_{a}\beta^{4-b}a}{2(1+\beta^2)a^*} \int _{\R^N\backslash B_\delta (0)} \frac{\bar{D}^{1+2\beta^2}_{a}}{|x|^{b}}\bar \xi_a,\\[4mm]
\end{array} \label{5.3:5}\end{equation}
where $ \bar{D}^{s-1}_{a}=\int^{1}_{0}\big[t\bar{u}_{2,a}+(1-t)\bar{u}_{1,a}\big]^{s-1}dt$ is defined by \eqref{dkp} and
\begin{equation}\arraycolsep=1.5pt\begin{array}{lll}
&&\displaystyle\frac{I_2-I_1}{\|\bar u_{2,a}-\bar u_{1,a}\|_{L^\infty}}\\[4mm]
&=&\displaystyle \frac{\alp^{b}_{a}\beta^{2-b}a}{a^*} \int_{\partial B_{\delta}(0)} \frac{(x\cdot\nu)}{|x|^{b}} \tilde{D}^{1+2\beta^2}_{a} \bar \xi_a  dS-\displaystyle\frac{\alp _a^2}{2} \int_{\partial B_{\delta}(0)} \big(\bar u_{2,a}+\bar u_{1,a}\big)\bar \xi_a V(x)(x\cdot\nu) dS\\[4mm]
&&+\displaystyle\frac{\mu_{2,a}\alp _a^2}{2} \int_{\partial B_{\delta}(0)}   \big(\bar u_{2,a}+\bar u_{1,a}\big)\bar \xi_a (x\cdot\nu) dS
+\displaystyle\frac{\big(\mu_{2,a}-\mu_{1,a}\big)\alp _a^2}{2\|\bar u_{2,a}-\bar u_{1,a}\|_{L^\infty}} \int_{\partial B_{\delta}(0)}   \bar u_{1,a}^2 (x\cdot\nu) dS.
\end{array} \label{5.3:10}
\end{equation}

We next estimate the right hand side of (\ref{5.3:5}). We first consider $\frac{I_2-I_1}{\|\bar u_{2,a}-\bar u_{1,a}\|_{L^\infty}}$. Using H\"{o}lder inequality, we derive from \eqref{5.5} and \eqref{na:A} that
\begin{equation}\arraycolsep=1.5pt\begin{array}{lll}
&&\displaystyle \frac{\alp^{b}_{a}\beta^{2-b}a}{a^*} \int_{\partial B_{\delta}(0)} \frac{(x\cdot\nu)}{|x|^{b}} \bar{D}^{1+2\beta^2}_{a} \bar \xi_a  dS-\displaystyle\frac{\alp _a^2}{2} \int_{\partial B_{\delta}(0)} \big(\bar u_{2,a}+\bar u_{1,a}\big)\bar \xi_a V(x)(x\cdot\nu) dS\\[4mm]
&&\quad+\displaystyle\frac{\mu_{2,a}\alp _a^2}{2} \int_{\partial B_{\delta}(0)}   \big(\bar u_{2,a}+\bar u_{1,a}\big)\bar \xi_a (x\cdot\nu) dS=o(e^{-\frac{C\delta}{\alp_a}}).
\end{array} \label{5.2}
\end{equation}
Moreover, we derive from (\ref{5.7}) below  that
\begin{equation}
\displaystyle\frac{\big|\mu_{2,a}-\mu_{1,a}\big|\alp _a^2}{\|\bar u_{2,a}-\bar u_{1,a}\|_{L^\infty}}
\le  \displaystyle\frac{2\beta^{N+4-b}a}{\alp^{N-b}_{a}\|Q\|^{2(1+\beta^2)}_{2}}\int_{\R^{N}} |x|^{-b}\bar{D}^{1+2\beta^2}_{a}|\bar \xi_a |
\le  C,
\label{5.2:9F}
\end{equation}
which combined with \eqref{5.2} yields that as $a\nearrow a^*$
\begin{equation}\label{5.3}
\displaystyle\frac{I_2-I_1}{\|\bar u_{2,a}-\bar u_{1,a}\|_{L^\infty}}=o(e^{-\frac{C\delta}{\alp_a}}).
\end{equation}
Moreover, if $\delta >0$ is small, we then deduce from \eqref{5.5} and \eqref{na:A} and \eqref{5.2:6} that
\begin{equation}\arraycolsep=1.5pt\begin{array}{lll}
&&\displaystyle \Big|\frac{2-N}{4}\alp _a^2\int_{\partial B_{\delta}(0)}\big(\nabla \bar u_{2,a}+\nabla \bar u_{1,a}\big)\nabla \bar \xi_a\Big|\\[4mm]
&&\displaystyle \leq C\alp _a^2\Big(\int_{\partial B_{\delta}(0)}|\nabla\bar{\xi}_{a}|^{2}\Big)^{\frac{1}{2}} \Big[\int_{\partial B_{\delta}(0)}\Big(|\nabla\bar{u}_{1,a}|^{2}\Big)^{\frac{1}{2}}+\int_{\partial B_{\delta}(0)}\Big(|\nabla\bar{u}_{2,a}|^{2}\Big)^{\frac{1}{2}}\Big]dS\\[4mm]
& &\displaystyle \leq C\alp _a^{\frac{N+2}{2}}e^{-\frac{C\delta}{\alp_a}},
\end{array} \label{1.1}\end{equation}
where $C>0$ is independent of $a\nearrow a^*$.  Similarly, we can estimate the rest terms in the right hand side of $T_a$.

We finally conclude from above  that
\[
T_a=o(e^{-\frac{C\delta}{\alp_a}}) \,\ \mbox{as} \,\ a\nearrow a^*,
\]
and thus we obtain from \eqref{5.3:4} that
\begin{equation}\arraycolsep=1.5pt\begin{array}{lll}
o(e^{-\frac{C\delta}{\alp_a}})=\displaystyle-\alp^{2}_{a}\int_{B_{\delta}(0)} \Big[V(x)+\frac{1}{2}[x\cdot \nabla V(x)]\Big](\bar u_{1,a}+\bar u_{2,a})\bar{\xi}_a.
\end{array} \label{5.3:12}
\end{equation}
Recall from the assumptions $(V_2)$ that if $\delta>0$ is small enough and $x\in B_\delta(0)$,  we have
\[
V(x)+\frac{1}{2}[x\cdot \nabla V(x)]=[1+o(1)]\big[h(x)+\frac{1}{2}[x\cdot \nabla h(x)]\big]=[1+o(1)]\frac{2+l}{2}h(x),
\]
where we used the fact that if $h(x)$ is a homogeneous function of degree $l>0$, then $x\cdot \nabla h(x)=lh(x)$ in the last equality. Note that $\tilde u_{i,a}(x)=\bar u_{i,a}\big(\frac{\alp_a}{\beta}\big)\to Q$ as $a\nearrow a^*$, we derive that as $a\nearrow a^*$
\[
\begin{aligned}
o(e^{-\frac{C\delta}{\alp_a}})=&\displaystyle-\alp^{2}_{a}\int_{B_{\delta}(0)} \Big[V(x)+\frac{1}{2}[x\cdot \nabla V(x)]\Big](\bar u_{1,a}+\bar u_{2,a})\bar{\xi}_adx\\
&=-\alp^{2}_{a}\int_{B_{\delta}(0)}[1+o(1)]\frac{2+l}{2}h(x) (\bar u_{1,a}(x)+\bar u_{2,a}(x))\bar{\xi}_adx\\
&=-\frac{\alp^{2+N+l}_{a}}{\beta^{N+l}}\int_{B_{\delta}(0)}[1+o(1)]\frac{2+l}{2}h(y) (\tilde u_{1,a}(y)+\tilde u_{2,a}(y))\xi_a(y)dy\\
&=-\frac{\alp^{2+N+l}_{a}}{\beta^{N+l}}\inte (2+l)h(y)Q(y)\xi_0(y)dy+o(\alp^{2+N+l}_{a}).
\end{aligned}
\]
Following \eqref{ks0}, we then conclude from above  that
\begin{equation}\arraycolsep=1.5pt\begin{array}{lll}
0&&=\displaystyle 2\int_{\R^{N}}h(x)Q(x)\xi_0\\[4mm]
&&=\displaystyle 2b_0\int_{\R^{N}}h(x)Q(x)\big(\frac{N}{2}Q+x\cdot\nabla Q\big)\\[4mm]
&&=\displaystyle 2b_0\Big[\frac{N}{2}\int_{\R^{N}}h(x)Q^2(x)+\frac{1}{2}\int_{\R^{N}}h(x)\big(x\cdot\nabla Q^2(x)\big)\Big]\\[4mm]
&&=\displaystyle 2b_0 \Big\{\frac{N}{2}\int_{\R^{N}}h(x)Q^2(x)- \frac{1}{2}\int_{\R^{N}} Q^2(x)\big[Nh(x)+x\cdot\nabla h(x)\big]\Big\}\\[4mm]
&&=\displaystyle -lb_0\int_{\R^{N}}h(x)Q^2(x),
\end{array} \label{5.4:5}
\end{equation}
which implies that $b_0=0$, i.e., $\xi_0\equiv0$ in $\R^N$.

\vskip 0.1truein
\noindent{\em  Step 3.}
$\xi_0\equiv 0$ cannot occur.

Finally, let $y_a$ be a point satisfying $|\xi_a(y_a)|=\|\xi_a\|_{L^\infty(\R^N)}=1$. Since $\tilde u_{i,a}$ ($i=1,2$) decays exponentially uniformly for $a\nearrow a^*$, applying the maximum principle to (\ref{3d:2}) yields that $|y_a|\le C$ uniformly in $a$.  Therefore, taking a subsequence if necessary, we assume $y_a\to y_0$ as $a\nearrow a^*$ and thus $|\xi_0(y_0)|=\lim_{a\nearrow a^{\ast}}|\xi_a(y_a)|=1 $, which however contradicts to the fact that $\xi_0 \equiv 0$ on $\R^N$. This completes the proof of Theorem \ref{unique}.
\qed

\appendix
\section{Appendix}
Firstly, we devoted to prove the equivalence between minimizers of (\ref{e(a)}) and the ground states of (\ref{F}). For convenience, we introduce some notations in advance. For any $a\in(0,a^*)$, where $a^*$ is defined by \eqref{a8}, the set of nontrivial weak solutions for (\ref{F}) is defined by
\begin{equation*}
S_{a,\mu}:=\{u\in \mathcal{H}\backslash\{0\}:\ \langle J'_{a,\mu}, \varphi\rangle=0 \ \mbox{for}\ \mbox{all}\  \varphi\in\mathcal{H}\},
\end{equation*}
where the energy functional $J_{a,\mu}$ is defined as
\begin{equation}\label{ja}
J_{a,\mu}:=\int_{\R^{N}}(|\nabla u|^{2}+\big(V(x)-\mu\big)|u(x)|^{2})dx-\frac{a}{1+\beta^2}\int_{\R^{N}}\frac{|u(x)|^{2+2\beta^2}}{|x|^{b}}dx.
\end{equation}
Furthermore, the set of ground states for (\ref{F}) is given by
\begin{equation}
G_{a,\mu}:=\{u\in  S_{a,\mu}:\ J_{a,\mu}(u)\leq J_{a,\mu}(v)\ \mbox{for}\ \mbox{all} \ v\in S_{a,\mu}\}.
\end{equation}
Moreover, we denote the set of minimizers for $e(a)$ as
\begin{equation}
\Lambda_a:=\{u_a\in\mathcal{H}(\R^N):\ u_a\ \mbox{is}\ \mbox{a}\ \mbox{minimizer}\ \mbox{of}\ e(a)\}.
\end{equation}
We have the following theorem.
\begin{thm}\label{A.3}
Suppose $V(x)$ satisfies $(V_1)$, then
\begin{enumerate}
\item For a.e. $a\in(0, a^*)$, all minimizers of $e(a)$ satisfy (\ref{F}) with a fixed Lagrange multiplier $\mu=\mu_a$.
\item For a.e. $a\in(0, a^*)$, $G_{a,\mu_a}=\Lambda_a$.
\end{enumerate}
\end{thm}
Since the proof of Theorem \ref{A.3} is very similar to that used in \cite[Theorem 1.1]{SX}, we omit it here for simplicity.

\begin{lem}\label{A.1}
	Suppose $w_a(x)\in H^1(\R^N)$, then $w^q_a(x)|x|^{-b}\in L^m(\R^N)$, where $m\in\Big(\frac{2N}{Nq+2b}, \frac{2^*N}{Nq+2^*b}\Big)$, $2^*=\frac{2N}{N-2}$ if $N\geq3$ and  $m\in\Big(\frac{2N}{Nq+2b}, \frac{N}{b}\Big)$ if $N=1,2$.
\end{lem}
\noindent\textbf{Proof.}
\begin{equation}\label{A1}
\begin{split}
\int_{\R^N}\big(w^q_a(x)|x|^{-b}\big)^mdx&=\int_{B_R(0)}\big(w^q_a(x)|x|^{-b}\big)^mdx+\int_{B^c_R(0)}\big(w^q_a(x)|x|^{-b}\big)^mdx\\
&:=J_1+J_2.
\end{split}
\end{equation}
By H\"{o}lder inequality,
\begin{equation*}
J_1\leq \Big(\int_{B_R(0)}w^{qms}_a(x)dx\Big)^{\frac{1}{s}}\Big(\int_{B_R(0)}|x|^{-bms'}dx\Big)^{\frac{1}{s'}},
\end{equation*}
\begin{equation*}
J_2\leq \Big(\int_{B^c_R(0)}(w^{qmt}_a(x)dx\Big)^{\frac{1}{t}}\Big(\int_{B^c_R(0)}|x|^{-bmt'}dx\Big)^{\frac{1}{t'}}.
\end{equation*}
For $N\geq 3$, $J_1<\infty$ if $m>0$ satisfies
\begin{equation}
\left\{
\begin{split}
&2\leq qms\leq 2^*,\\
&0<bs'm<N,\\
&s>1,\,s'=\frac{s}{s-1},
\end{split}
\right.
\end{equation}
which are equivalent to $m>0$ satisfies
\begin{equation}\label{m1}
\left\{
\begin{split}
&\frac{2}{qs}\leq m\leq\frac{2^*}{qs}, \\
&0< m<\frac{N(s-1)}{bs},\\
&s>1.
\end{split}
\right.
\end{equation}
Now let
\begin{equation}
m_1(s):=\frac{2}{qs},\quad m_2(s):=\frac{2^*}{qs},\quad  m_3(s):=\frac{N(s-1)}{bs},\,\ \hbox{where $s>1$,}
\end{equation}
then from (\ref{m1}), $m$ is located in the intersection part between the middle of $m_1(s)$, $m_2(s)$ and the below of $m_3(s)$, which implies $m\in \big(0,\frac{2^*N}{Nq+2^*b}\big)$.

On the other hand, $J_2<\infty$ if $m>0$ satisfies
\begin{equation}
\left\{
\begin{split}
&2\leq qmt< 2^*\\
&bt'm>N,\\
&t>1,\,t'=\frac{t}{t-1},
\end{split}
\right.
\end{equation}
which are equivalent to $m>0$ satisfies
\begin{equation}\label{t1}
\left\{
\begin{split}
&\frac{2}{qt}\leq m<\frac{2^*}{qt}, \\
&m>\frac{N(t-1)}{bt},\\
&t>1.
\end{split}
\right.
\end{equation}
Let
\begin{equation}
m_1(t):=\frac{2}{qt},\quad m_2(t):=\frac{2^*}{qt},\quad m_3(t):=\frac{N(t-1)}{bt},\,\ \hbox{where $t>1$,}
\end{equation}
then by (\ref{t1}), $m$ is located in the intersection part between the middle of $m_1(t)$, $m_2(t)$ and the above of $m_3(t)$, which implies $m\in \big(\frac{2N}{Nq+2b}, \frac{2^*}{q}\big)$. 

Therefore, we conclude from the above two results that $w^q_a(x)|x|^{-b}\in L^m(\R^N)$ if $m\in \big(0,\frac{2^*N}{Nq+2^*b}\big)\cap \big(\frac{2N}{Nq+2b}, \frac{2^*}{q}\big)$, i.e., $m\in\Big(\frac{2N}{Nq+2b}, \frac{2^*N}{Nq+2^*b}\Big)$.
When $N=1,2$, note that $2^*=\infty$ in this case. By the same method as the above, we can obtain that $w^q_a(x)|x|^{-b}\in L^m(\R^N)$ if $m\in \big(\frac{2N}{Nq+2b},\frac{N}{b}\big)$ in this case. The proof of Lemma \ref{A.1} is therefore completed.\qed

\vskip 0.1truein
\noindent{\bf Proof of the claim \eqref{5.2:6}.} Denote
\begin{equation}\label{dkp}
\arraycolsep=1.5pt\begin{array}{lll}
\bar{D}^{s-1}_{a}(x):&&=\displaystyle\frac{\bar{u}_{2,a}^{s}(x)-\bar{u}_{1,a}^{s}(x)}{s(\bar{u}_{2,a}-\bar{u}_{1,a})}\\[4mm]
&&\displaystyle=\frac{\int^{1}_{0}\frac{d}{dt}[t\bar{u}_{2,a}+(1-t)\bar{u}_{1,a}]^{s}dt}{s(\bar{u}_{2,a}-\bar{u}_{1,a})}=\int^{1}_{0}\big[t\bar{u}_{2,a}+(1-t)\bar{u}_{1,a}\big]^{s-1}dt.
\end{array}
\end{equation}
then we obtain from \eqref{xy} \eqref{5.6} that $\bar \xi_a$ satisfies the following equation
\begin{equation}
-\alp ^2_a\Delta \bar \xi_a +\bar C_{a}(x)\bar \xi_a =\bar g_a(x)\quad \text{in\,\, $\R^N$},
\label{5.2:1}
\end{equation}
where the coefficient $\bar C_a(x)$ satisfies
\begin{equation}
\bar C_a(x):=-\mu _{1,a}\alp ^2_a-\frac{(1+2\beta^2)\alp^{b}_{a}\beta^{2-b}a}{a^*}|x|^{-b}\bar{D}^{2\beta^2}_{a}+ \alp ^2_a  V(x),
\label{5.2:2}
\end{equation}
and the nonhomogeneous term $\bar g_a(x)$ satisfies
\begin{equation}
\arraycolsep=1.5pt\begin{array}{lll}
\bar g_a(x)&:=\displaystyle \frac{\alp ^2_a\bar u_{2,a}(\mu _{2,a}-\mu _{1,a})}{\|\bar u_{2,a}-\bar u_{1,a}\|_{L^\infty(\R^N)}}\\[4mm]
&=-\displaystyle  \frac{2\beta^{N+4-b}a\bar u_{2,a}}{\|Q\|^{2(1+\beta^2)}_{2}\alp^{N-b}_a}\int_{\R^{N}} \bar \xi_a |x|^{-b}\bar{D}^{1+2\beta^2}_{a}dx,
\end{array}
\label{5.2:3}
\end{equation}
where used the following fact in the last equality, 
\begin{equation}\label{5.7}
\arraycolsep=1.5pt\begin{array}{lll}
\mu _{2,a}-\mu _{1,a}&=-\displaystyle\frac{\beta^2a}{1+\beta^2}\int_{\R^{N}}|x|^{-b}(u^{2(1+\beta^2)}_{2,a}-u^{2(1+\beta^2)}_{1,a})
\\[4mm]
&=-\displaystyle\frac{\beta^2a}{1+\beta^2}\frac{\beta^{N(1+\beta^2)}}{\alp_{a}^{N(1+\beta^2)}\|Q\|_2^{2(1+\beta^2)}}\int_{\R^{N}}|x|^{-b}(\bar{u}^{2(1+\beta^2)}_{2,a}-\bar{u}^{2(1+\beta^2)}_{1,a})
\\[4mm]
&=-\displaystyle\frac{2\beta^{N+4-b}a\|\bar u_{2,a}-\bar u_{1,a}\|_{L^\infty(\R^N)}}{\|Q\|_{2}^{2(1+\beta^2)}\alp_a^{N(1+\beta^2)}}\int_{\R^{N}} \bar \xi_a |x|^{-b}\bar{D}^{1+2\beta^2}_{a}dx.\end{array}
\end{equation}
Multiplying (\ref{5.2:1}) by $\bar \xi_a$ and integrating over $\R^N$, we obtain that
\[\arraycolsep=1.5pt\begin{array}{lll}
&&\displaystyle \alp ^2_a\int_{\R^{N}} |\nabla \bar \xi_a|^2 -\mu_{1,a}\alp ^2 _a\int_{\R^{N}}  |\bar\xi_a|^2+\alp ^2_a\int_{\R^{N}} V(x)|\bar\xi_a|^2 \\[4mm]
&=&\displaystyle \frac{(1+\beta^2)\alp^{b}_{a}\beta^{2-b}a}{a^*}\int_{\R^{N}} |x|^{-b}\bar{D}^{2\beta^2}_{a}(x)|\bar\xi_a|^2\displaystyle-\frac{a\beta^{N+4-b}}{\alp^{N-b}_{a}\|Q\|^{2(1+\beta^2)}_{2}}\int_{\R^{N}} \bar u_{2,a}\bar\xi_a\int_{\R^{N}}\bar\xi_a |x|^{-b}\bar D^{1+2\beta^2}_{a}(x)\\[4mm]
&\leq &\displaystyle \frac{(1+2\beta^2)\alp^{b}_{a}\beta^{2-b}a}{a^*}\int_{\R^{N}}|x|^{-b}\bar{D}^{2\beta^2}_{a}(x)
\displaystyle+\frac{2a\beta^{N+4-b}}{\alp^{N-b}_{a}\|Q\|^{2(1+\beta^2)}_{2}}\int_{\R^{N}} \bar u_{2,a}\int_{\R^{N}}|x|^{-b}\bar{D}^{1+2\beta^2}_{a}(x)\\[4mm]
&\le & C\alp ^N_a\,\ \mbox{as} \,\ a\nearrow a^*,
\end{array}\]
due to the fact that $|\bar\xi_a|$ and $\bar u_{i,a}\big(\frac{\alp _a}{\beta}x\big)$ are bounded uniformly for $a$, and $\bar u_{i,a}\big(\frac{\alp _a}{\beta}x\big)$ satisfies \eqref{na:A}, where $i=1,\, 2$. Recall from \eqref{uyp} that $\alp_a^2\mu_{i,a}\to-\beta^2$ as $a\nearrow a^*$, the above estimate further implies that there exists a constant $C_1>0$ such that
\begin{equation}
I:=\alp ^2_a\int_{\R^{N}} |\nabla \bar \xi_a|^2+\frac{\beta ^2}{2} \int_{\R^{N}} |\bar\xi_a|^2+ \alp ^2_a\int_{\R^{N}} V(x)|\bar\xi_a|^2<C_1\alp ^N_a \quad \text{as}\ \, a\nearrow a^*.
\label{5.2:5}
\end{equation}
Following  \cite[Lemma 4.5]{Cao}, we then conclude that for any $x_0\in\R^N$, there exist a small constant $\delta >0$ and $C_2>0$  such that
\[
\int_{\partial B_\delta (x_0)} \Big[ \alp ^2_a |\nabla \bar \xi_a|^2+ \frac{\beta ^2}{2}  |\bar\xi_a|^2+ \alp ^2_a  V(x)|\bar\xi_a|^2\Big]dS\le C_2I\le C_1C_2\alp ^N_a\,\ \mbox{as} \,\ a\nearrow a^*,
\]
and the claim (\ref{5.2:6}) is thus proved. \qed

\begin{lem}\label{A.2}
Let $\mathcal{L}=-\Delta +1-(1+2\beta^2\big)|x|^{-b}Q^{2\beta^2}$, where $Q(x)$ is the unqiue positive radially symmetric solution of $$-\Delta Q+Q-Q^{2\beta^2}Q|x|^{-b}=0\,\ \hbox{in $\R^N$},$$ 
then $\mathcal{L}\big(\frac{N}{2}Q+x\nabla Q\big)=-2Q$.
\end{lem}
\textbf{Proof.}
Since $-\Delta Q+Q-Q^{2\beta^2}Q|x|^{-b}=0$, taking the partial derivative of both sides of the equation, we have
\begin{equation}
-\Delta (\partial_iQ)+\partial_iQ-\big(1+2\beta^2\big)Q^{2\beta^2}\partial_iQ|x|^{-b}+bQ^{2\beta^2}Q|x|^{-(b+2)}x_i=0,
\end{equation}
and hence $\mathcal{L}\big(\partial_iQ\big)=-bQ^{2\beta^2}Q|x|^{-(b+2)}x_i$.
Moreover, by simple calculation, we get
\begin{equation}
\begin{split}
\mathcal{L}\big(x_i\partial_iQ\big)&=-\Delta (x_i\partial_iQ)+x_i\partial_iQ-(1+2\beta^2\big)|x|^{-b}Q^{2\beta^2}x_i\partial_iQ\\
&=x_i\mathcal{L}\big(\partial_iQ\big)-2\partial_i\big(\partial_iQ\big).
\end{split}
\end{equation}
We then deduce from the above facts that
\begin{equation}
\begin{split}
\mathcal{L}\big(x\cdot\nabla Q\big)&=-b|x|^{-b}Q^{2\beta^2}Q-2\Delta Q\\
&=-b|x|^{-b}Q^{2\beta^2}Q-2\big(Q-|x|^{-b}Q^{2\beta^2}Q\big)\\
&=(2-b)|x|^{-b}Q^{2\beta^2}Q-2Q.
\end{split}
\end{equation}
On the other hand,
\begin{equation}
\begin{split}
\mathcal{L}\big(Q\big)&=-\Delta Q+Q-(1+2\beta^2\big)|x|^{-b}Q^{2\beta^2}Q\\
&=-\Delta Q+Q-|x|^{-b}Q^{2\beta^2}Q-2\beta^2|x|^{-b}Q^{2\beta^2}Q\\
&=-2\beta^2|x|^{-b}Q^{2\beta^2}Q\\
&=-\frac{4-2b}{N}|x|^{-b}Q^{2\beta^2}Q,
\end{split}
\end{equation}
where we have used the facts $-\Delta Q+Q-Q^{2\beta^2}Q|x|^{-b}=0$ and $\beta^2=\frac{2-b}{N}$. Therefore we conclude that $\mathcal{L}\big(\frac{N}{2}Q+x\cdot\nabla Q\big)=-2Q$ and the proof of Lemma \ref{A.2} is completed.\qed

\end{document}